
\def\LaTeX{%
  \let\Begin\begin
  \let\End\end
  \let\salta\relax
  \let\finqui\relax
  \let\futuro\relax}

\def\UK{\def\our{our}\let\sz s}
\def\USA{\def\our{or}\let\sz z}

\UK 



\LaTeX

\USA


\salta

\documentclass[twoside,12pt]{article}
\setlength{\textheight}{24cm}
\setlength{\textwidth}{16cm}
\setlength{\oddsidemargin}{2mm}
\setlength{\evensidemargin}{2mm}
\setlength{\topmargin}{-15mm}
\parskip2mm


\usepackage[usenames,dvipsnames]{color}
\usepackage{amsmath}
\usepackage{amsthm}
\usepackage{amssymb,bbm}
\usepackage[mathcal]{euscript}

\usepackage{cite}
\usepackage{hyperref}
\usepackage{enumitem}

\usepackage[ulem=normalem,draft]{changes}
%
%

%
 
\definecolor{ciclamino}{rgb}{0.5,0,0.5}
\definecolor{blu}{rgb}{0,0,0.7}
\definecolor{rosso}{rgb}{0.85,0,0}

\def\juerg #1{{\color{red}#1}}
\def\an #1{{\color{magenta}#1}}
\def\gianni #1{{\color{Green}#1}}
\def\pier #1{{\color{blue}#1}}
 
\def\rev #1{{\color{rosso}#1}}

\def\rev #1{#1}
\def\gianni #1{#1}
\def\juerg #1{#1}
\def\an #1{#1}
\def\pier #1{#1}




\bibliographystyle{plain}


%

\finqui

\def\Beq{\Begin{equation}}
\def\Eeq{\End{equation}}
\def\Bsist{\Begin{eqnarray}}
\def\Esist{\End{eqnarray}}

\def\Bthm{\Begin{theorem}}
\def\Ethm{\End{theorem}}
\def\Blem{\Begin{lemma}}
\def\Elem{\End{lemma}}

\def\Bcor{\Begin{corollary}}
\def\Ecor{\End{corollary}}
\def\Brem{\Begin{remark}\rm}
\def\Erem{\End{remark}}

\def\Bdim{\Begin{proof}}
\def\Edim{\End{proof}}
\def\Bcenter{\Begin{center}}
\def\Ecenter{\End{center}}
\let\non\nonumber




\def\step #1 \par{\medskip\noindent{\bf #1.}\quad}
\def\jstep #1: \par {\vspace{2mm}\noindent\underline{\sc #1 :}\par\nobreak\vspace{1mm}\noindent}

\def\aand{\quad\hbox{and}\quad}
\def\Lip{Lip\-schitz}
\def\Holder{H\"older}

\def\lhs{left-hand side}
\def\rhs{right-hand side}
\def\sfw{straightforward}




\def\multibold #1{\def\arg{#1}%
  \ifx\arg\pto \let\next\relax
  \else
  \def\next{\expandafter
    \def\csname #1#1\endcsname{{\boldsymbol #1}}%
    \multibold}%
  \fi \next}

\def\pto{.}

\def\multical #1{\def\arg{#1}%
  \ifx\arg\pto \let\next\relax
  \else
  \def\next{\expandafter
    \def\csname cal#1\endcsname{{\cal #1}}%
    \multical}%
  \fi \next}

\def\multigrass #1{\def\arg{#1}%
  \ifx\arg\pto \let\next\relax
  \else
  \def\next{\expandafter
    \def\csname grass#1\endcsname{{\mathbb #1}}%
    \multigrass}%
  \fi \next}


\def\multimathop #1 {\def\arg{#1}%
  \ifx\arg\pto \let\next\relax
  \else
  \def\next{\expandafter
    \def\csname #1\endcsname{\mathop{\rm #1}\nolimits}%
    \multimathop}%
  \fi \next}

\multibold
qweryuiopasdfghjklzxcvbnmQWERTYUIOPASDFGHJKLZXCVBNM.  

\multical
QWERTYUIOPASDFGHJKLZXCVBNM.

\multigrass
QWERTYUIOPASDFGHJKLZXCVBNM.

\multimathop
diag dist div dom mean meas sign supp .

\def\Span{\mathop{\rm span}\nolimits}


\def\accorpa #1#2{\eqref{#1}--\eqref{#2}}
\def\Accorpa #1#2 #3 {\gdef #1{\eqref{#2}--\eqref{#3}}%
  \wlog{}\wlog{\string #1 -> #2 - #3}\wlog{}}


\def\separa{\noalign{\allowbreak}}

\def\somma #1#2#3{\sum_{#1=#2}^{#3}}

\def\graffe #1{\mathopen\{#1\mathclose\}}
\def\<#1>{\mathopen\langle #1\mathclose\rangle}
\def\norma #1{\mathopen \| #1\mathclose \|}

\def\aeQ{\checkmmode{a.e.\ in~$Q$}}
\def\aet{\checkmmode{a.e.\ in~$(0,T)$}}
\def\aat{\checkmmode{for a.a.\ $t\in(0,T)$}}

\let\hat\widehat
\def\cpto{\,\cdot\,}

\def\Beta{\widehat\beta}
\def\Pi{\widehat\pi}
\def\Betaeps{\Beta_\eps}
\def\betaeps{\beta_\eps}

\def\phieps{\phi_\eps}
\def\mueps{\mu_\eps}
\def\weps{w_\eps}

\def\iot {\int_0^t}
\def\ioT {\int_0^T}
\def\intQt{\int_{Q_t}}
\def\intQ{\int_Q}
\def\iO{\int_\Omega}

\def\dt{\partial_t}
\def\dtt{\partial_t^2}
\def\dn{\partial_{\nn}}

\def\0{{\boldsymbol {0} }}

\let\emb\hookrightarrow

\def\checkmmode #1{\relax\ifmmode\hbox{#1}\else{#1}\fi}


\let\erre\grassR
\let\enne\grassN




\def\genspazio #1#2#3#4#5{#1^{#2}(#5,#4;#3)}
\def\spazio #1#2#3{\genspazio {#1}{#2}{#3}T0}

\def\L {\spazio L}
\def\H {\spazio H}
\def\W {\spazio W}

\def\C #1#2{C^{#1}([0,T];#2)}


\def\Lx #1{L^{#1}(\Omega)}
\def\Hx #1{H^{#1}(\Omega)}
\def\Wx #1{W^{#1}(\Omega)}

\def\LQ #1{L^{#1}(Q)}

\def\Luno{\Lx 1}
\def\Ldue{\Lx 2}

\def\Huno{\Hx 1}
\def\Hdue{\Hx 2}



\let\theta\vartheta

\let\eps\varepsilon

\let\phi\varphi

\let\TeXchi\chi                         
\newbox\chibox
\setbox0 \hbox{\mathsurround0pt $\TeXchi$}
\setbox\chibox \hbox{\raise\dp0 \box 0 }
\def\chi{\copy\chibox}



\def\ej{e_j}
\def\ei{e_i}
\def\lambdaj{\lambda_j}
\def\Vn{V_n}
\def\phin{\phi_n}
\def\mun{\mu_n}
\def\wn{w_n}
\def\phinj{\phi_{nj}}
\def\munj{\mu_{nj}}
\def\wnj{w_{nj}}
\def\vn{v_n}

\def\hphin{\hat {\boldsymbol{\phi}}_n}
\def\hmun{ \hat {\boldsymbol{\mu}}_n}
\def\hwn{ \hat {\boldsymbol{w}}_n}

\def\ubar{\overline u}
\def\vbar{\overline v}

\def\psibar{\overline\psi}
\def\phibar{\overline\phi}

\def\phinbar{\overline\phin}
\def\munbar{\overline\mun}
\def\fbar{\overline f}
\def\phizbar{\overline\phiz}

\def\Pn{\grassP_n}

\def\normaV #1{\norma{#1}_V}

\def\normaVp #1{\norma{#1}_*}

\def\phiz{\phi_0}
\def\wz{w_0}
\def\wu{w_1}

\def\Vp{{V^*}}

\def\soluz{(\phi,\mu,\xi,w)}
\def\soluzeps{(\phieps,\mueps,\weps)}
\def\soluzn{(\phin,\mun,\wn)}

\def\CO{C_\Omega}
\def\kuno {{\kappa_1}}
\def\kdue {{\kappa_2}}
\def\rmin{r_*}
\def\rmax{r^*}


\usepackage{amsmath}
\DeclareFontFamily{U}{mathc}{}
\DeclareFontShape{U}{mathc}{m}{it}%
{<->s*[1.03] mathc10}{}

\DeclareMathAlphabet{\mathscr}{U}{mathc}{m}{it}

\Begin{document}


%
\title{
	\juerg{On a Cahn--Hilliard system with source term and thermal memory}
}
\author{}
\date{}
\maketitle
\Bcenter
\vskip-1.5cm
{\large\sc Pierluigi Colli$^{(1)}$}\\
{\normalsize e-mail: {\tt pierluigi.colli@unipv.it}}\\[0.25cm]
{\large\sc Gianni Gilardi $^{(1)}$}\\
{\normalsize e-mail: {\tt gianni.gilardi@unipv.it}}\\[0.25cm]
{{\large\sc Andrea Signori$^{(2)}$}\\
{\normalsize e-mail: {\tt andrea.signori@polimi.it}}}\\[0.25cm]
{\large\sc J\"urgen Sprekels$^{(3)}$}\\
{\normalsize e-mail: {\tt juergen.sprekels@wias-berlin.de}}\\[.5cm]
$^{(1)}$
{\small Dipartimento di Matematica ``F. Casorati'', Universit\`a di Pavia}\\
{\small and Research Associate at the IMATI -- C.N.R. Pavia}\\
{\small via Ferrata 5, I-27100 Pavia, Italy}\\[.3cm] 
{
$^{(2)}$
{\small Dipartimento di Matematica, Politecnico di Milano}\\
{\small via E. Bonardi 9, I-20133 Milano, Italy}
}
\\[.3cm] 
$^{(3)}$
{\small Department of Mathematics}\\
{\small Humboldt-Universit\"at zu Berlin}\\
{\small Unter den Linden 6, D-10099 Berlin, Germany}\\
{\small and}\\
{\small Weierstrass Institute for Applied Analysis and Stochastics}\\
{\small Mohrenstrasse 39, D-10117 Berlin, Germany}\\[10mm]
\date{}

\Ecenter
\Begin{abstract}
\noindent 
A nonisothermal phase field system \juerg{of Cahn--Hilliard type} is introduced and analyzed mathematically. \juerg{The system constitutes an 
extension of the classical Caginalp model for nonisothermal phase transitions with a conserved order parameter. It couples a Cahn--Hilliard type equation with source term for the order parameter with the universal balance law of internal energy. In place of the standard 
Fourier form, the constitutive law of the heat flux is assumed in the form given by the theory developed by Green and Naghdi, which 
accounts for a possible thermal memory of the evolution. This has the consequence that the balance law of internal energy becomes
a second-order in time equation for the {\it thermal displacement} or {\it freezing index}, that is, a primitive with respect to time of
the temperature. Another particular feature of our system is the presence of the source term in the equation for the order parameter, which
entails additional mathematical difficulties because the mass conservation of the order parameter is lost.} 
We provide several mathematical results under general assumptions on the source term and the double-well nonlinearity \juerg{governing the evolution}:
existence and continuous dependence results are shown for weak and strong solutions to the corresponding initial-boundary value problem.


\End{abstract}

\pagestyle{myheadings}
\newcommand\testopari{\sc Colli -- Gilardi -- Signori -- Sprekels}
\newcommand\testodispari{\sc {Cahn--Hilliard system with source term and thermal memory}}
\markboth{\testopari}{\testodispari}
%

\section{Introduction}
\label{INTRO}
\setcounter{equation}{0}

A common assumption in phase segregation processes of binary mixtures is to postulate that the mixture under investigation undergoes \juerg{the phase 
separation at a constant temperature. However, in numerous applications the evolution does not take place under isothermal conditions.}
The first contribution aiming at including temperature effects in the theory of phase separation is due to Caginalp \cite{Cag1, Cag2, Cag3}.
\juerg{It} was motivated by the Stefan problem for the evolution of \juerg{the interface in a solid-liquid} phase transition and \juerg{in a}
Hele--Shaw type flow between two fluids with different viscosities.

Another typical assumption in the context of the Cahn--Hilliard equation is the mass conservation property that arises as a direct 
consequence of the standard no-flux \juerg{boundary} condition prescribed for the chemical potential associated with the phase field variable.
\juerg{While this condition is} very natural for the engineering applications that Cahn and Hilliard had in mind originally (see \cite{CH}), 
the recent employment of the Cahn--Hilliard equation to describe other \juerg{phenomena driven by phase segregation} 
demands the incorporation of an external source term $S$
\juerg{in the model that reflects} the fact that the system may not be isolated and the loss or production of mass is possible.
Without claiming to be exhaustive, let us mention that numerous liquid-liquid phase segregation problems arise in \pier{cell biology}~\cite{Bio}
and in tumor growth models~\cite{GLSS}. 
For this reason, we also included the presence of a source term in our investigation.
 
The \juerg{standard isothermal Cahn--Hilliard system has been extensively studied} in the past decades: see, e.g., \cite{CHbook} 
and the references therein. 
On the other hand, the mathematical understanding of nonisothermal Cahn--Hilliard systems is, 
\juerg{thirty years after the seminal works by Alt and Pawlow (see \cite{AP,AP2} and, in particular, \cite{AP1}) and twenty years 
after the groundbreaking work \cite{Gaj} by Gajewski for the nonlocal case, still \juerg{far from being complete.}
Before presenting our system, let us discuss some recent literature.}
Concerning some analytic results of the aforementioned system by Caginalp, we mention the related contibutions \cite{CC12,CC13,Mir}.
Next, employing micro-force balance theory, Miranville and Schimperna proposed a further derivation in \cite{MS}, and the well-posedness of a related system has been addressed in \cite{SM}.
Moreover, we point out the recent contribution \cite{Chun} by De Anna et al., where two new thermodynamically consistent models related to 
nonisothermal Cahn--Hilliard systems have been derived. 
Finally, we refer to \cite{Gal1,Gal2} for some mathematical results on a relaxed version of the above systems endowed with dynamic boundary conditions.

Motivated by the aforementioned remarks, we aim at analyzing a nonisothermal Cahn--Hilliard \pier{type} system with source term in this paper.
To this end, let $\Omega\subset \erre^3$ be the spatial domain where the evolution takes place, and $T>0$ a given final time.
\juerg{We then consider the following initial-boundary value problem}:
\begin{alignat}{2}
  & \dt\phi - \Delta\mu + \gamma \phi
  = f  
  \quad && \text{in $Q:=\Omega \times (0,T)$},
  \label{Iprima}
  \\
  & \mu
  = - \Delta\phi + F'(\phi) + a - b \dt w
  \quad && \text{in $Q$},
  \label{Iseconda}
  \\
  & \dtt w - \Delta(\kuno \dt w + \kdue w) + \lambda \dt\phi
  = g
  \quad && \text{in $Q$},
  \label{Iterza}
  \\
  & \rev{\dn\phi=0,
  \quad
  \dn\mu
  = 0,
  \quad 
  \dn(\kuno \dt w + \kdue w)
  = 0 }
  \quad && \text{on $\Sigma:=\partial\Omega \times (0,T)$},
  \label{Ibc}
  \\
  & \phi\rev{\vert_{t=0}}= \phiz, 
  \quad 
  w\rev{\vert_{t=0}} = \wz,
  \quad 
  \dt w\rev{\vert_{t=0}}= \wu 
  \quad && \text{in $\Omega.$}
  \label{Icauchy}
\end{alignat}
\Accorpa\Iequaz Iprima Iterza
\Accorpa\Ipbl Iprima Icauchy
\rev{In the above system several  positive physical constants carrying physical dimensions
 have, for convenience, been taken equal to unity, while keeping their physical dimensions.
This  has no bearing on the mathematical analysis on which we focus here. The unknowns have the following meaning:}
{$\phi$ is a normalized difference between the volume fractions
of pure phases in the binary mixture (the dimensionless {\em order parameter} of the phase transformation, which should attain its values in 
the interval $[-1,1]$), 
\rev{and $\mu$ corresponds to the \emph{chemical potential}; the unknown $w$ stands for
the so-called \emph{thermal displacement}, which in the mathematical literature of free boundary problems is also termed \emph{freezing index} or \emph{Baiocchi transform}. It is directly connected}
to the temperature $\theta$ (which in the case of the Caginalp model is actually 
a temperature difference) through the relation}
\begin{align}
	\label{thermal_disp}
	w (\cdot , t)  = w_0 + \iot \theta(\cdot, s) \,ds, \quad t \in[0,T].
\end{align}
Moreover, $\kuno$ and $\kdue$ in \eqref{Iterza} stand for prescribed positive coefficients related to the heat flux\pier{;}
$\gamma$ is a positive physical constant related to the intensity of the mass absorption/production of the source\pier{, where the source term in \eqref{Iprima} is $f- \gamma \varphi$ as explained below;}
$\lambda$ stands for the latent heat of the phase transformation\pier{;} $a, b$ are physical constants\pier{;} $g$ is  a distributed heat source.
Besides, the symbol $\dn$ represents the outward normal derivative on $\Gamma:= \partial\Omega$, while   $\phi_0,w_0,$ and $w_1$ indicate some given initial values. 
Finally, $F'$ stands for the (generalized) derivative of a double-well shaped nonlinearity. 
Prototypical and important examples for $F$ 
are the so-called {\em classical regular potential} and the {\em logarithmic double-well potential},
which are the functions given~by
\Bsist
  && F_{reg}(r) := \frac 14 \, (r^2-1)^2 \,,
  \quad r \in \erre ,
  \label{regpot}
  \\
  &&   
	\label{logpot}
  F_{log}(r) := 
  \left\{
    \begin{array}{ll}
      (1+r)\ln (1+r)+(1-r)\ln (1-r) - c_1 r^2 
      & \quad \hbox{if $|r|\leq1$},
      \\
      +\infty
      & \quad \hbox{otherwise},
	\end{array}
  \right.
\Esist
with the convention $0\ln(0):=\lim_{r\searrow0}r\ln(r)=0$.
In \eqref{logpot}, $c_1>1$ so that $F_{log}$ is nonconvex.
Another example is {the {\em double obstacle potential\/}, where, with} $c_2>0$,
\Beq
  F_{2obs}(r) := - c_2 r^2 
  \quad \hbox{if $|r|\leq1$}
  \aand
  F_{2obs}(r) := +\infty
  \quad \hbox{if $|r|>1$}.
  \label{obspot}
\Eeq
Singular potentials like \eqref{logpot} and \eqref{obspot} are \juerg{difficult to handle from the mathematical viewpoint, but   have the great advantage that if a solution exists, then it automatically inherits
the property of being physically meaningful}, that is, $\phi \in [-1,1]$. In general, this cannot be guaranteed for regular potentials 
like the quartic \eqref{regpot}, which, in this sense, provides just an approximation of the more physical choices.
In cases like \eqref{obspot}, one has to split $F$ into a nondifferentiable convex part $\hat\beta$
(the~indicator function of $[-1,1]$ in the present example) and a smooth (usually quadratic) perturbation $\hat\pi$.
Accordingly, the second equation \eqref{Iseconda} has then to be understood as the differential inclusion
\begin{align*}
	\mu \in - \Delta\phi + \partial \hat\beta(\phi) + \hat \pi'(\phi) + a - b \dt w,
\end{align*}
or, equivalently, with the help of a selection $\xi$, as the identity
\begin{align*}
	\mu= - \Delta\phi + \xi + \pi(\phi) + a - b \dt w
	\quad \text{with} \quad 
	\xi \in \partial\hat\beta(\phi).
\end{align*}

\juerg{
The above system is a formal extension of the Cahn--Hilliard system introduced by Caginalp in \cite{Cag2}
(see also the derivation in \cite[Ex.~4.4.2, (4.44), (4.46)]{BS}); it corresponds to the Allen--Cahn counterpart analyzed 
recently in \cite{CSS3}. The main differences between our system and the one originally introduced in \cite{Cag2}
are the following:} 

\juerg{%
\noindent
$\bullet$ \,\,In \cite{Cag2}, we have $\,a=\lambda\,$ (the specific latent heat).\\[2mm]
$\bullet$ \,\,In \cite{Cag2}, the heat flux is assumed in the standard Fourier form $\mathbf{q}=-\kappa_1\nabla\theta$, while we 
we follow the works by Green and Naghdi~\cite{GN91,GN92,GN93} and Podio-Guidugli~\cite{PG09} \rev{(see also \cite{SFPG})} and postulate that
\begin{equation}\label{flux}
\mathbf q=-\kappa_1 \nabla (\dt w )- \kappa_2 \nabla w \quad\mbox{where $\kappa_i>0$, $i=1,2$.}
\end{equation} 
Note that this assumption accounts for a possible previous
thermal history of the phenomenon. We also observe that the no-flux condition $\mathbf q\cdot\mathbf n=0$ then gives rise to the third
boundary condition in \eqref{Ibc}.\\[2mm]
\noindent
$\bullet$ \,\,The third -- and main -- difference is that \eqref{Iprima}--\eqref{Iseconda} comprises a Cahn--Hilliard system 
with  a source term \,$S:=f - \gamma\phi$, which is independent of temperature.}

\juerg{The presence of $S$ radically changes the behavior of the Cahn--Hilliard equation since} the mass conservation property is no longer fulfilled.
In fact, due to the no-flux boundary condition for $\mu$ in \eqref{Ibc}, a formal consideration, i.e., testing \eqref{Iprima} by $1/|\Omega|$, readily reveals that
the mass balance law of the order parameter $\phi$ is ruled by
\begin{align*}
	 \frac d{dt} \Big(\frac 1 {|\Omega|}\iO \phi(t) \Big)
	= 
	\frac 1 {|\Omega|} \iO {S(t)}
	\quad \aat.
\end{align*}
In this direction, we highlight that, especially when working with singular potential like \eqref{obspot}, the control of the mean value of $\phi$ plays a crucial role in the mathematical analysis of Cahn--Hilliard-type systems.
Besides, in the case when $f$ is a positive constant such that $f \in (-\gamma ,\gamma)$, the equations \eqref{Iprima}--\eqref{Iseconda} correspond to the so-called Cahn--Hilliard--Oono system, see, e.g., \cite{GGM} and \cite{CGRS4}.

Finally, let us mention that the differential structure with respect to $w$ in equation \eqref{Iterza} is sometimes also referred to as the {\it strongly damped wave equation} (see, e.g., \cite{Pata} and the references therein).

Let us conclude this section by presenting an outline of the paper.
In the following section, we state the main results and list the corresponding assumptions.
Then, from Section~\ref{UNIQUENESS} onward, we start proving the mentioned results. In particular, Section~\ref{UNIQUENESS} is devoted to showing some continuous dependence results enjoyed by the system \Ipbl.
In Section~\ref{APPROXIMATION}, we then introduce and solve a preparatory approximating problem that will allow us to prove in Section~\ref{EXISTENCE} the existence of weak solutions, as well as some regularity results.


\section{Statement of the problem and \an{main} results}
\label{STATEMENT}
\setcounter{equation}{0}

Throughout the paper, $\Omega$ \an{indicates} a bounded and connected open subset of~$\erre^3$
(the lower dimensional cases can be treated in the same way)
\juerg{with smooth boundary $\Gamma:=\partial\Omega$.}
In the following, $|\Omega|$~and $\nn$ denote the Lebesgue measure of $\Omega$
and the outward unit normal vector field on~$\Gamma$, respectively.
Given a final time $T>0$, we~set
\begin{align}
	Q_t := \Omega \times (0,t)
	\quad \hbox{for $t\in(0,T]$}
	\aand
	Q := Q_T \,.
	\label{defQt}
\end{align}
Given a Banach space $X$, we denote its norm by $\norma\cpto_X$,
with the exceptions of \an{$L^p$} spaces on $\Omega$ and~$Q$,
whose norms are denoted by $\norma\cpto_p$ for $1\leq p\pier{{}\leq\infty{}}$
(if~no confusion can arise), and of the space $H$ introduced below.
For brevity, we use the same symbol for the norm in a space and in any power \juerg{thereof}.
\an{
Furthermore, for Banach spaces $X$ and $Y$, we notice that  the linear space $X\cap Y$ becomes a Banach space when equipped
with its natural graph norm 
\begin{align*}
	\norma{v}_{X\cap Y}:= \norma{v}_X + \norma{v}_Y, 
	\quad 
	v \in X \cap Y.
\end{align*}
}
Then, we introduce the shorthands
\Beq
  H := \Ldue , \quad
  V := \Huno,
  \aand
  W := \graffe{v\in\Hx2: \ \dn v=0 \an{\,\text{ on $\Gamma$}}},
  \label{defspazi}
\Eeq
and endow these spaces with their natural norms.
For simplicity, we write $\norma\cpto$ instead of~$\norma\cpto_H$.
Moreover, we denote by $(\cpto,\cpto)$ and $\<\cpto,\cpto>$ 
the standard inner product of $H$ and the duality pairing between the dual space $\Vp$ of $V$ and $V$ itself, respectively.
We identify $H$ with a subspace of $\Vp$ in the usual way, i.e.,
in order that 
\an{
\begin{align*}
	\text{$\<u,v>=(u,v)$ 
	\quad 
	for every $u\in H$ and $v\in V$.}
\end{align*}
}
This makes $(V,H,\Vp)$ a Hilbert triplet.
We notice that all of the embeddings
\Beq
  W \emb V \emb H \emb V^*
  \non
\Eeq
are dense and compact.
Next, we define the generalized mean value $\vbar$ of a generic element $v\in\Vp$ by setting
\Beq
  \vbar := \frac 1{|\Omega|} \, \< v , 1 >\,,
  \label{defmean}
\Eeq
where we have written $1$ for the constant function that takes the \an{value $1$ in $\Omega$}.
It is clear that $\vbar$ reduces to the usual mean value if $v\in H$.
The same notation $\vbar$ is \an{employed also} if $v$ is a time-dependent function.

Let us come to the structural assumptions we make for our analysis.
First, 
\Beq
  \hbox{$\gamma$, $a$, $b$, $\kuno$, $\kdue$ and $\lambda$ are positive constants}.
  \label{hpconst}
\Eeq
Next, in order to allow \an{for} general \an{double-well} potentials \an{in \eqref{Iseconda}}, we assume~that
\Bsist
  && F : \erre \to (0,+\infty]
  \quad \hbox{admits the decomposition} \quad
  F = \Beta + \Pi ,
  \quad \hbox{where}
  \label{hpF}
  \\[1mm]
  && \Beta : \erre \to [0,+\infty]
  \quad \hbox{is convex\an{,} l.s.c.\an{,} and fulfills} \quad
  \Beta(0)= 0 ,
  \label{hpBeta}
  \\[1mm]
  && \Pi \in C^1(\erre),
  \quad \hbox{and its derivative is \Lip\ continuous}.
  \label{hpPi}
\Esist
Moreover, we set
\Beq
  \beta := \partial\Beta
  \aand
  \pi := \Pi',
  \label{defbetapi}
\Eeq
\Accorpa\Hpstruttura hpconst defbetapi
\an{where $\partial$ denotes the subdifferential operator,} and notice that $\beta:\erre\to2^{\erre}$ is maximal monotone with corresponding domain $D(\beta)$
and that $0\in\beta(0)$.
We observe that all of the examples \accorpa{regpot}{obspot} of potentials \an{introduced before do}
satisfy the conditions required \an{above}.
\an{
Of course, in the case of 
nonregular potentials like the double obstacle \eqref{obspot}, the second equation \eqref{Iseconda} has to be intended as the differential inclusion
\begin{align*}
	\mu \in - \Delta\phi + \beta(\phi) + \pi(\phi) + a - b \dt w,
\end{align*}
or, equivalently, with the help of a selection $
	\xi  \in \beta(\phi) \,\,\aeQ$, as the identity
\begin{align*}
	\mu= - \Delta\phi + \xi + \pi(\phi) + a - b \dt w.
\end{align*}
}

As for the data, we assume \pier{that}
\Bsist
  && f \in \LQ\infty
  \aand
  g \in \L2H ,
  \label{hpfg}
  \\
  && \phiz \in W , \quad
  \wz \in V 
  \aand
  \wu \in H \,.
  \label{hpz}
\Esist
However, we also need some compatibility conditions 
between the data $f$ and $\phiz$ and the domain of~$\beta$.
\an{\juerg{These} are in fact already expected, as we are dealing with possible singular potentials
and a mass source.
In particular, let us repeat that the contribution \gianni{$S:=f - \gamma\phi$} in \eqref{Iprima} 
plays the role of \juerg{a} (phase-dependent) mass source/sink in the model. 
Indeed, by formally testing \eqref{Iprima} by $1/|\Omega|$, and using \eqref{Ibc}, 
we infer that the mass balance law of the system reads
\begin{align*}
	 \frac d{dt} \Big(\frac 1 {|\Omega|}\iO \phi(t) \Big)
	= 
	\frac 1 {|\Omega|} \iO \gianni{S(t)}
	\quad \gianni\aat.
\end{align*}
Therefore, it is \juerg{natural to expect} to have some compatibility conditions 
between the structure of the source term $S$, thus on the constant $\gamma$ and the function $f$,
and the possibly singular potential $\beta$.}
Namely, setting
\Beq
  \rho := \frac {\norma f_\infty} \gamma\,,
  \label{defrho}
\Eeq
and noting that $\phiz\in\juerg{C^0(\overline\Omega)}$ by \eqref{hpz}, we 
\juerg{require~that all of the quatities}
\Bsist
  && \min_{x\in\overline\Omega}\phiz(x) , \ 
  \max_{x\in\overline\Omega}\phiz(x), \
  - \rho - (\phizbar)^- \,, \ \rho + (\phizbar)^+
  \non
  \\
  && \quad \hbox{belong to the interior of $D(\beta)$},
  \label{compatibility}
\Esist
\Accorpa\Hpdati hpfg compatibility
where $(\cpto)^-$ and $(\cpto)^+$ denote the negative and positive part functions\an{, respectively}.

\Brem
\label{Remdati}
The assumptions on $f$ and $\phiz$ \juerg{can} be weakened \juerg{slightly}.
However, \juerg{in} doing so, we \juerg{would have to} replace \accorpa{defrho}{compatibility}
by more complicated \an{compatibility} conditions.
Moreover, \juerg{when} regularizing our problem as we \juerg{are going to} do in \an{the forthcoming} Section~\ref{APPROXIMATION},
we \juerg{would have to} regularize $\phiz$ as well.
This would lead to estimates depending on the regularization parameter, \an{so that} further uniform estimate\an{s} 
\juerg{had to} be performed.
\Erem

At this point, we can \an{rigorously state} \juerg{our notion} of \an{(weak)} solution to the \an{aforementioned} problem \juerg{under study}.
A \an{weak} solution \an{to the system \Ipbl} is a quadruplet $\soluz$ enjoying the regularity \an{properties}
\Bsist
  && \phi \in \H1\Vp \cap \L\infty V \cap \L2W,
  \label{regphi}
  \\
  && \mu \in \L2V,
  \label{regmu}
  \\
  \separa
  && \xi \in \L2H,
  \label{regxi}
  \\
  && w \in \H2\Vp \cap \W{1,\infty}H \cap \H1V,
  \label{regw}
\Esist
\Accorpa\Regsoluz regphi regw
and \juerg{satisfying}
\Bsist
  && \< \dt\phi , v >
  + \iO \nabla\mu \cdot \nabla v
  + \gamma \iO \phi v
  = \iO f v
  \non
  \\
  && \quad \hbox{for every $v\in V$ and \aet}\,,
  \label{prima}
  \\
  \separa
  && \mu
  = - \Delta\phi + \xi + \pi(\phi) + a - b \dt w
  \aand 
  \xi \in \beta(\phi)
  \quad \aeQ\,,
  \label{seconda}
  \\
  \separa
  && \< \dtt w , v >
  + \iO \nabla( \kuno \dt w + \kdue w) \cdot \nabla v
  + \lambda \iO \dt\phi \, v
  = \iO g v
  \non
  \\
  && \quad \hbox{for every $v\in V$ and \aet}\,,
  \label{terza}
  \\
  && \phi\rev{\vert_{t=0}} = \phiz, 
  \quad 
  w\rev{\vert_{t=0}} = \wz,
  \aand
  \dt w\rev{\vert_{t=0}} = \wu \,.
  \label{cauchy}
\Esist
\Accorpa\Pbl prima cauchy	 

The present paper is devoted to \juerg{the study of the} well-posedness of the above problem
and of the regularity of its solutions.
\juerg{Our first result}  is an existence theorem.

\Bthm
\label{Existence}
Assume \Hpstruttura\ on the structure of the system and \Hpdati\ on the data.
Then, problem \Pbl\ has at least one solution $\soluz$ 
satisfying \Regsoluz\ \juerg{and}
\Beq
  \phi \in \L2{\Wx{2,6}}
  \aand
  \xi \in \L2{\Lx6}\,,
  \label{piureg}
\Eeq
\juerg{as well as} the estimate
\Bsist
  && \norma\phi_{\H1\Vp \cap \L\infty V \cap \L2{\Wx{2,6}} }
  +\norma\mu_{\L2V}
  + \norma\xi_{\L2{\Lx6}}
  \non
  \\
  &&  {}
  + \norma w_{\H2\Vp \cap \W{1,\infty}H \cap \H1V}
  \leq K_1\,,
  \label{stima}
\Esist
with a \an{positive} constant $K_1$ that depends only on the structure of the system,
$\Omega$, $T$, and upper bounds for the norms of the data 
and the quantities related to assumptions \Hpdati.
\Ethm

Uniqueness cannot be expected, in general, as it usually occurs in Cahn\an{--}Hilliard type problems \an{with nonregular potentials}.
However, we have the result stated below, 
which ensures continuous dependence on $f$ and $g$
for the components $\phi$ and $w$ of every solution
with fixed initial data.
In the statement, \juerg{we use} the \an{following} notation \pier{for convolution products with $1$}:
\Beq
  (1*v)(t) := \iot v(s) \, ds
  \quad \hbox{for $v\in\L1H$ \,and\, $t\in [0,T]$} .
  \label{convoluz}
\Eeq

\Bthm
\label{Contdep}
Under the assumptions \Hpstruttura\ on the structure of the system and 
\an{\accorpa{hpz}{compatibility}} on the initial data,
let $f_i$ and~$g_i$, $i=1,2$, satisfy \eqref{hpfg}, and let
$(\phi_i,\mu_i,\xi_i,w_i)$ be any two corresponding solutions \pier{of problem \Pbl} \gianni{with the regularity \Regsoluz}.
Then, the estimate
\Bsist
  && \norma{\phi_1-\phi_2}_{\L\infty\Vp\cap\L2V}
  + \norma{w_1-w_2}_{\H1H\cap\L\infty V}
  \non
  \\
  && \leq K_2 \bigl(
    \norma{f_1-f_2}_{\L2\Vp\cap\LQ1}
    + \norma{f_1-f_2}_{\LQ1}^{1/2}
    + \norma{1*(g_1-g_2)}_{\L2H}
  \bigr)
  \qquad
  \label{contdep}
\Esist
holds true with a \an{positive} constant $K_2$ that depends only on the structure of the system,
$\Omega$, $T$, and an upper bound for the norms of $\xi_1$ and $\xi_2$ in~$\LQ1$.
\Ethm

Partial uniqueness in general and full uniqueness if $\beta$ is single valued trivially follow,
as stated below.

\Bcor
\label{Uniqueness}
Assume \Hpstruttura\ on the structure of the system and \Hpdati\ on the data.
Then, the components $\phi$ and $w$ of any solution  \juerg{in the sense of Theorem 2.2} are uniquely determined.
Furthermore, if $\beta$ is single valued, \juerg{then}
even the components $\mu$ and $\xi$ are uniquely determined 
and the solution is unique.
\Ecor

Under proper regularity assumption on $\beta$ and on the data, there exists a more regular solution.
We notice that all of the examples \accorpa{regpot}{obspot} of potentials \an{still}
satisfy the \an{stronger} condition\an{s required below}.

\Bthm
\label{Regularity}
In addition to the assumptions of Theorem~\ref{Existence}, \juerg{let the following conditions be fulfilled:}
\begin{align}
  &	\juerg{\hbox{\pier{the} restriction of $\beta$ to the interior of $D(\beta)$ is a single-valued $C^1$-function,}}
  \label{hpbetareg}
  \\[1mm]
  & f \in \pier{{}\H{1}\Vp{}}
  , \quad
  \phiz \in \Hx3\pier,   
  \aand
  \wu \in V \an{.}
  \label{hpdatireg}
\end{align}
Then the problem \Pbl\ admits at least one solution $\soluz$ that enjoys the further regularity
\begin{align}
  & \phi \in \H1V \cap \L\infty{\Wx{2,6}}, \quad
  \mu \in \L\infty V ,
  \non
  \\
  &
  \xi \in \L\infty{\Lx6},
  \quad
  \gianni{w \in \H2H \cap \W{1,\infty}V,}
  \label{regularity}
\end{align}
and satisfies the estimate
\Bsist
  && \norma\phi_{\H1V\cap\L\infty{\Wx{2,6}}}
  + \norma\mu_{\L\infty V}   + \norma\xi_{\L\infty{\Lx6}}  \non
  \\
  && \quad {}
  + \gianni{\norma w_{\H2H\cap\W{1,\infty}V}}
  \,\leq \,K_3\,,
  \label{stimareg}
\Esist
with a \gianni{positive constant} $K_3$ that depends only on the structure of the system,
$\Omega$, $T$, and upper bounds of the norms of the data  
and the quantities related to assumptions \Hpdati\ and \eqref{hpdatireg}.
\Ethm

\gianni{%
\Brem
Notice that \eqref{terza} says that
$u:=\kuno\dt w+\kdue w$ satisfies
\Beq
  \iO \nabla u \cdot \nabla v
  = \iO h v
  \quad \hbox{for every $v\in V$ and \aet},
  \non
\Eeq
where $h:=g-\lambda\dt\phi-\dtt w$.
In particular, if $\soluz$ is a solution \juerg{in the sense of} Theorem~\ref{Regularity}, \juerg{then}
$h$ belongs to $\L2H$, and the elliptic regularity theory yields that
\Beq
  u \in \L2W
  \aand
  \norma u_{\L2W}
  \leq \CO \bigl( \norma u_{\L2V} + \norma h_{\L2H} \bigr),
  \non
\Eeq
where $\CO$ depends only on~$\Omega$.
Thus, the same norm can be estimated by a constant \juerg{which is} proportional to~$K_3$.
By solving $\kuno\dt w+\kdue w=u$ for~$w$, 
we obtain that $w\in\H1\Hdue$ or $w\in\H1W$, provided that $\wz\in\Hdue$ or $\wz\in W$, respectively.
\Erem
}%

In \an{the next sections, when} proving our results, we \juerg{widely use  \Holder's inequality,} 
as well as the Young, Poincar\'e, Sobolev and compactness \an{inequalities} \juerg{recalled below:}
\Bsist
  && ab \leq \delta a^2 + \frac 1{4\delta} \, b^2
  \quad \hbox{for every $a,b\in\erre$ and $\delta>0$}.
  \label{young}
  \\[2mm]
  \separa
  && \normaV v
  \leq \CO \, \bigl( \norma{\nabla v} + |\vbar| \bigr)
  \quad \hbox{for every $v\in V$}.
  \label{poincare}
  \\[2mm]
  \separa
  && \norma v_p
  \leq \CO \, \normaV v
  \quad \hbox{for every $v\in V$ and $p\in[1,6]$}.
  \label{sobolev}
  \\[2mm]
  \separa
  && \norma v_p 
  \leq \delta \, \norma{\nabla v} + C_{\Omega,p,\delta} \, \normaVp v
  \quad \hbox{for every $v\in V$, $p\in[1,6)$ and $\delta>0$}.
  \label{compact}
\Esist
\juerg{Here,} $\CO$ is a constant that depends only on~$\Omega$,
while $C_{\Omega,p,\delta}$ depends on $p$ and~$\delta$, in addition.
Moreover, the symbol $\normaVp\cpto$ denotes the norm in $\Vp$ 
defined \an{by} the forthcoming formula~\eqref{normaVp}.
The last two inequalities are related to the (three-dimensional) embedding
$V\emb\Lx p$ which holds for $p\in[1,6]$ and is compact if $p<6$.

Finally, we take advantage of a \juerg{tool that is commonly used} \an{in} the study of problems related to the Cahn--Hilliard \an{type} equations:
consider, for $\psi\in\Vp$, the problem of finding
\Beq
  u \in V
  \quad \hbox{such that} \quad
  \iO \nabla u \cdot \nabla v
  = \< \psi , v >
  \quad \hbox{for every $v\in V$}.
  \label{neumann}
\Eeq
\juerg{Obviously,} if $\psi\in H$, \juerg{then this problem is just} the usual homogeneous Neumann problem
for the \an{Poisson equation} $-\Delta u=\psi$.
Now, since $\Omega$ is connected, for $\psi\in\Vp$, \eqref{neumann} is solvable if and only if $\psi$ has zero mean value.
Moreover, if this condition is satisfied, \juerg{then}
\an{there exists a unique solution} \an{possessing} zero mean value.
This \an{entails} that the operator
\Bsist
  && \calN: \dom(\calN) := \graffe{\psi\in\Vp:\ \psibar=0} \to \graffe{u\in V:\ \ubar=0},
  \ \hbox{ given by the rule}
  \non
  \\
  && \quad \psi \mapsto \hbox{the unique solution $u$ to \eqref{neumann} satisfying $\ubar=0$},
  \qquad\qquad
  \label{defN}
\Esist
is well defined and \juerg{yields} an isomorphism between the above spaces.
\an{Besides, it} follows that the formula
\Beq
  \normaVp\psi^2 := \norma{\nabla\calN(\psi-\psibar)}^2 + |\psibar|^2
  \quad \hbox{for $\psi\in\Vp$}
  \label{normaVp}
\Eeq
defines a Hilbert norm in $\Vp$ that is equivalent to the standard dual norm.
From the above definitions one trivially derives~that
\Bsist
  && \iO \nabla\calN\psi \cdot \nabla v
  = \< \psi , v >
  \quad \hbox{for every $\psi\in\dom(\calN)$ and $v\in V$},
  \label{dadefN}
  \\[2mm]
  && \< \psi , \calN\zeta >
  = \< \zeta , \calN\psi >
  \quad \hbox{for every $\psi,\zeta\in\dom(\calN)$},
  \label{simmN}
  \\[2mm]
  \separa
  && \< \psi , \calN\psi > 
  = \iO |\nabla\calN\psi|^2
  = \normaVp\psi^2
  \quad \hbox{for every $\psi\in\dom(\calN)$}.
  \label{danormaVp}
\Esist
Moreover, it turns out that
\Beq
  \iot \< \dt v(s) , \calN v(s) > \, ds
  = \iot \< v(s) , \calN(\dt v(s)) > \, ds
  = \frac 12 \, \normaVp{v(t)}^2
  - \frac 12 \, \normaVp{v(0)}^2
  \label{propN} 
\Eeq
for every $t\in[0,T]$ and every $v\in\H1\Vp$ satisfying $\vbar=0$ \aet.
\Accorpa\PropN defN propN

We conclude this section by stating a general rule concerning the constants 
that appear in the estimates to be performed in the following.
The small-case symbol $\,c\,$ stands for a generic constant
whose actual value may change from line to line, and even within the same line,
and depends only on~$\Omega$, the shape of the nonlinearities,
and the constants and the norms of the functions involved in the assumptions of the statements.
In particular, the values of $\,c\,$ do not depend on the parameters $\,\eps>0\,$ and $\,n\in\enne\,$ 
that will be introduced in the next sections.
A~small-case symbol with a subscript like $c_\delta$ 
indicates that the constant may depend on the parameter~$\delta$, in addition.
On the contrary, we mark precise constants that we can refer~to
by using different symbols
(see, e.g., \eqref{sobolev}).


\section{\gianni{Continuous dependence}}
\label{UNIQUENESS}
\setcounter{equation}{0}

This section is devoted to the proof of \an{T}heorem~\ref{Contdep}.
\gianni{Let $f_i$ and~$g_i$, $i=1,2$, satisfy \eqref{hpfg}, and let
$(\phi_i,\mu_i,\xi_i,w_i)$ be any two corresponding solutions as in the statement.
We set, for convenience, $\phi:=\phi_1-\phi_2$,} 
\juerg{and define $\mu$, $\xi$, $w$, $f$ and~$g$ analogously}.
We first make some preliminary observations.
\an{Recalling \gianni{\eqref{regphi} (see \eqref{defspazi} for the definition of~$W$) 
and testing~{\eqref{prima}} by~$ 1/{|\Omega|}$}, we find}~that
\Beq
  \frac d{dt} \, \phibar(t) 
  + \gamma \, \phibar(t)
  = \fbar(t)
  \quad \aat.
  \label{meanphi}
\Eeq
Then, \juerg{on the one hand}, by multiplying this equality by $\iO v$,
we deduce~that
\Beq
  \iO \dt\phibar \, v 
  + \gamma \iO \phibar \, v
  = \iO \fbar v
  \quad \hbox{for every $v\in V$ and \aet} \,.
  \label{meanprima}
\Eeq
On the other hand, by (formally) multiplying \eqref{meanphi} by $\pier{\sign(\phibar)}$, 
where $\sign:\erre\to\erre$ is the sign function defined by $\sign( r):=r/|r|$ if $r\not=0$ and $\sign(0)=0$,
we infer~that
\Beq
   |\phibar(t)|
  + \gamma \iot |\phibar(s)| \, ds
  \leq \iot |\fbar(s)| \, ds\,,
  \non
\Eeq
\pier{whence}
\Beq
  \sup_{t\in(0,T)}|\phibar(t)| 
  \leq \ioT |\fbar(s)| \, ds
  \leq \frac 1{|\Omega|} \, \norma f_{\LQ1} \,.
  \label{stimaphibar}
\Eeq

We now start the proof of the theorem.
We use the properties \PropN\ of the operator~$\calN$ 
and recall the \pier{notation}~\eqref{convoluz} \pier{for convolution products with $1$}.
We write \an{equation} \eqref{prima} for both solutions and take the difference,
\juerg{obtaining} an equality from which we subtract~\eqref{meanprima} \juerg{to arrive at the identity}
\Beq
  \< \dt(\phi-\phibar) , v >
  + \iO \nabla\mu \cdot \nabla v
  + \gamma \iO (\phi-\phibar) v
  = \iO (f-\fbar) v
  \non
\Eeq
for every $v\in V$ and \aet.
\juerg{Since} $(\phi-\phibar)(t)$ has zero mean value for every $t\in[0,T]$,
we \an{are allowed to} test the above equation by $\calN(\phi-\phibar)$\an{. Integration} with respect to time \an{then leads to}, for every $t\in[0,T]$,
\Bsist
  && \frac 12 \, \normaVp{(\phi-\phibar)(t)}^2
  + \intQt \mu (\phi-\phibar)
  + \gamma \iot \normaVp{\phi-\phibar}^2
  = \intQt (f-\fbar) \, \calN(\phi-\phibar).
  \qquad
  \label{diffprima}
\Esist
Next, we write \eqref{seconda} for both solutions, multiply the difference by $-(\phi-\phibar)$,
and integrate over~$Q_t$\juerg{, finding} that
\Bsist
  && \intQt |\nabla\phi|^2
  + \intQt \xi \, \phi
  \an{{}-\intQt \mu (\phi-\phibar)}
  \an{{}- b \intQt \dt w \, \phi}
  \non
  \\
  && 
  = \intQt \xi \, \phibar
  - \intQt \bigl( \pi(\phi_1)-\juerg{\pi}(\phi_2) \bigr) (\phi-\phibar)
  \an{{}- b \intQt \dt w \, \phibar} \,.
  \label{diffseconda}
\Esist
Finally, we write \eqref{terza} for both solutions and take the convolution with~$1$.
Then, we test the difference of the \an{corresponding equalities} by $(b/\lambda)\dt w$ \an{to obtain}~that
\Bsist
  && \frac b\lambda \intQt |\dt w|^2 
  + \frac {b\kuno}{2\lambda} \iO |\nabla w(t)|^2
  \an{{}+{}} b \intQt \phi \, \dt w
  \non
  \\
  && 
  = - \frac {b\kdue} \lambda \intQt \nabla(1*w) \cdot \nabla\dt w
  + \an{\frac b\lambda}\intQt (1*g) \dt w \,.
  \label{diffterza}
\Esist
At this point, we add \accorpa{diffprima}{diffterza} to each other
and notice \an{that} some cancellations \an{occur}.
\juerg{Moreover, $\beta$ is monotone, and thus all of} the remaining terms on the \lhs\ are nonnegative.
We treat those on the \rhs\ \juerg{individually}.
First, we have~that
\Bsist
  && \intQt (f-\fbar) \, \calN(\phi-\phibar)
  \leq c \iot \normaVp{(f-\fbar)(s)} \, \normaV{\calN(\phi-\phibar)(s)} \, ds
  \non
  \\
  && \leq c \iot \normaVp{(f-\fbar)(s)} \, \normaVp{(\phi-\phibar)(s)} \, ds
  \non
  \\
  && \leq \iot \normaVp{(\phi-\phibar)(s)}^2 \, ds
  + c \, \norma{f-\fbar}_{\L2\Vp}^2 
  \leq \iot \normaVp{(\phi-\phibar)(s)}^2 \, ds
  + c \, \norma f_{\L2\Vp}^2 \,,
  \non
\Esist
where we have used the trivial inequalities
$\,\normaVp\vbar\leq c\,|\vbar|\leq c\,\normaVp v$,
which hold for every $v\in\Vp$.

Next, \juerg{we fix a constant $M$ such that $\norma{\xi_i}_{\LQ1}\leq M$ for $i=1,2$.
Then,} recalling \eqref{stimaphibar},
we have~that
\Beq
  \intQt \xi \, \phibar
  \leq \intQt \bigl( |\xi_1|+|\xi_2| \bigr) |\phibar|
  \leq 2M \! \sup_{s\in(0,t)}|\phibar(s)|
  \leq \frac{2M}{|\Omega|} \, \norma f_{\LQ1} \,.
  \non
\Eeq
\juerg{Also, in view of} the \Lip\ continuity of~$\pi$,
the obvious inequality $\norma\vbar\leq\norma v$ for $v\in H$,
and the compactness inequality \eqref{compact}, we find~that
\Beq
  - \intQt \bigl( \pi(\phi_1)-\juerg{\pi}(\phi_2) \bigr) (\phi-\phibar)
  \leq c \intQt |\phi|^2
  \leq \frac 12 \intQt |\nabla\phi|^2
  + c \iot \normaVp{\phi(s)}^2 \, ds \,.
  \non
\Eeq
Moreover\an{, by Young's \juerg{inequality} and arguing as above,} \pier{we have that}
\Beq
  - b \intQt \dt w \, \phibar
  \leq \frac b{4\lambda} \intQt |\dt w|^2 
  + c \intQt |\phibar|^2
  \leq \frac b{4\lambda} \intQt |\dt w|^2 
  + c \, \norma f_{\LQ1}^2\,,
  \non
\Eeq
on account of \eqref{stimaphibar}.
We deal with the next integral \an{using} integration by \juerg{parts} \an{to infer that}
\Bsist
  && - \frac {b\kdue} \lambda \intQt \nabla(1*w) \cdot \nabla\dt w
  = \frac {b\kdue} \lambda \intQt |\nabla w|^2
  - \frac {b\kdue} \lambda \iO \nabla(1*w)(t) \cdot \nabla w(t)
  \non
  \\
  && \leq  \an{c} \intQt |\nabla w|^2
  + \frac {b\kuno}{4\lambda} \iO |\nabla w(t)|^2
  + c \iO |\nabla(1*w)(t)|^2 \,.
  \non
\Esist
In addition, \pier{it is clear that}
\Beq
  \iO |\nabla(1*w)(t)|^2
  = \iO \Bigl| \iot \nabla w(s) \, ds \Bigr|^2
  \leq \iO t \iot |\nabla w(s)|^2 \, ds
  \leq T \intQt |\nabla w|^2 \,.
  \non
\Eeq
Finally, \pier{we note that}
\Beq
  \frac b\lambda \intQt (1*g) \dt w
  \leq \frac b{4\lambda} \intQt |\dt w|^2 
  + c \int_{\an{Q_t}} |1*g|^2 \,.
  \non
\Eeq
\an{Upon} collecting \an{\accorpa{diffprima}{diffterza}} and \juerg{the \an{inequalities} shown above}, we obtain~that
\Bsist
  && \frac 12 \, \normaVp{(\phi-\phibar)(t)}^2
  + \frac 12 \intQt |\nabla\phi|^2
  + \frac b{2\lambda} \intQt |\dt w|^2 
  + \frac {b\kuno}{4\lambda} \iO |\nabla w(t)|^2
  \non
  \\
  && \leq c \bigl(
  \norma f_{\L2\Vp\cap\LQ1}^2 
  + \norma f_{\LQ1}
  + \norma{1*g}_{\L2H}^2
  \bigr) 
  \non
  \\
  && 
  + c \Bigl(
    \iot \normaVp{\phi(s)}^2 \, ds
    + \intQt |\nabla w|^2 
  \Bigr)\,, 
  \non
\Esist
where $c$ has the dependence required for the constant $K_2$ \juerg{in the statement of the theorem}.
On the other hand, \eqref{stimaphibar} implies~that
\Bsist
  \normaVp{\phi(t)}
  \leq \normaVp{(\phi-\phibar)(t)}
  + c \, |\phibar(t)|
  \leq \normaVp{(\phi-\phibar)(t)}
  + c \, \norma f_{\LQ1}
  \quad \an{\aat}.
  \non
\Esist
By combining this with the previous inequality,
we are in \juerg{a} position to apply Gronwall's lemma and \pier{obtain} the desired estimate~\eqref{contdep}, \juerg{which concludes the proof}.


\section{Approximation}
\label{APPROXIMATION}
\setcounter{equation}{0}

In this section, we introduce and solve a proper approximating problem depending on the parameter~$\eps\in(0,1)$.
First of all, we replace the functional $\Beta$ and the \an{maximal monotone} graph $\beta$ by their Moreau--Yosida regularizations
$\Betaeps$ and~$\betaeps$, respectively
(see, e.g., \cite[pp.~28 and~39]{Brezis}).
We recall~that
\Bsist
  && 0 \leq \Betaeps(r)
  = \int_0^r \betaeps(s) \, ds
  \leq \Beta(r)
  \quad \hbox{for every $r\in\erre$}\,,
  \label{disugBetaeps}
  \\
  && \hbox{$\betaeps$ is monotone and \Lip\ continuous with $\betaeps(0)=0$}\,,
  \label{monbetaeps}
  \\[1mm]
  && |\betaeps(r)| \leq |\beta^\circ(r)|
  \quad \hbox{for every $r\in D(\beta)$}\,,
  \label{disugbetaeps}
\Esist
where $\beta^\circ(r)$ \an{denotes} the element of \an{the section} $\beta(r)$ having minimum modulus.
The approximating problem to be considered consists in finding a triplet $\soluzeps$
satisfying the regularity properties
\Bsist
  && \phieps \in \H1\Vp \cap \L\infty V \cap \L2{W\cap\Wx{2,6}}\,,
  \label{regphieps}
  \\
  && \mueps \in \L2V\,,
  \label{regmueps}
  \\
  && \weps \in \H2\Vp \cap \W{1,\infty}H \cap \H1V\,,
  \label{regweps}
\Esist
\Accorpa\Regsoluzeps regphieps regweps
and solving the following \gianni{system of variational identities or equations and initial conditions}:
\Bsist
  && \< \dt\phieps , v >
  + \iO \nabla\mueps \cdot \nabla v
  + \gamma \iO \phieps v
  = \iO f v
  \non
  \\
  && \quad \hbox{for every $v\in V$ and \aet}\,,
  \label{primaeps}
  \\
  \separa
  && \mueps
  = - \Delta\phieps + \betaeps(\phieps) + \pi(\phieps) + a- b \dt\weps
  \quad \aeQ\,,
  \label{secondaeps}
  \\
  \separa
  && \< \dtt\weps , v >
  + \iO \nabla( \kuno \dt\weps + \kdue\weps) \cdot \nabla v
  + \lambda \iO \dt\phieps \, v
  = \iO g v
  \non
  \\
  && \quad \hbox{for every $v\in V$ and \aet}\,,
  \label{terzaeps}
  \\
  && \phieps\rev{\vert_{t=0}} = \phiz, 
  \quad 
  \weps\rev{\vert_{t=0}} = \wz
  \aand
  \dt\weps\rev{\vert_{t=0}} = \wu \,.
  \label{cauchyeps}
\Esist
\Accorpa\Pbleps primaeps cauchyeps

\an{We remark that here we obviously do not need to consider any selection $\xi$ as $\betaeps$ is regular and single valued.}
Here is our basic result.

\Bthm
\label{Existenceeps}
Let the assumptions of Theorem~\ref{Existence} be in force.
Then problem \Pbleps\ has\an{, for every $\eps \in (0,1)$,} a unique solution \an{$\soluzeps$}
satisfying the regularity properties \an{expressed in} \Regsoluzeps.
\Ethm

The rest of this section is devoted to the proof of the above theorem.
Clearly, uniqueness is a consequence of Theorem~\ref{Contdep},
since $\betaeps$ satisfies all the assumptions \juerg{postulated} for $\beta$ \an{in \eqref{hpF}--\eqref{hpPi}},
and it is single valued, in addition \an{(cf. Corollary \ref{Uniqueness})}.

To prove the existence of a solution, we start from a Faedo--Galerkin scheme.
To this end, we introduce the nondecreasing \an{(ordered)} sequence
$\graffe{\lambdaj}$ of eigenvalues
and the corresponding complete orthonormal \juerg{sequence} $\graffe{\ej}$ of eigenfunctions
of the eigenvalue problem for the Laplace operator 
with homogeneous Neumann boundary conditions.
Namely, we have~that
\Bsist
  && - \Delta\ej = \lambdaj\ej
  \quad \hbox{in $\Omega$}
  \aand
  \dn\ej = 0
  \quad \hbox{on $\Gamma$}
  \quad \hbox{for $j=1,2,\dots$},
  \label{eigenj}
  \\
  && \iO \ei \ej 
  = \delta_{ij} 
  \quad \hbox{for every $i$ and $j$},
  \label{orthonormal}
\Esist
with the standard Kronecker symbol\an{s} $\delta_{ij}$.
Moreover, we~set
\Beq
  \Vn := \Span\graffe{\ej:\ 1\leq j\leq n}
  \quad \hbox{for $n=1,2,\dots$}
  \label{defspazin}
\Eeq
and recall that the union of these spaces is dense in both $V$ and~$H$.
Notice that all of the eigenfunctions are smooth since $\Omega$ is smooth.
\pier{Furthermore, as} $\Omega$ is connected, we have that $\lambda_1=0<\lambda_2$, 
and  $V_1$ is the subspace of \an{constant} functions.

The discrete problem consists \an{then} in finding a triplet $\soluzn$ of functions satisfying
\Beq
  \phin \in \H1\Vn \,, \quad
  \mun \in \L2\Vn
  \aand
  \wn \in \H2\Vn\,,
  \label{regsoluzn}
\Eeq
and solving the  \an{discrete} problem
\Bsist
  && \iO \dt\phin \, v
  + \iO \nabla\mun \cdot \nabla v
  + \gamma \iO \phin v
  = \iO f v
  \non
  \\
  && \quad \hbox{for every $v\in\Vn$ and \aet}\,,
  \label{priman}
  \\
  \separa
  && \iO \mun v
  = \iO \nabla\phin \cdot \nabla v
  + \iO \betaeps(\phin) v 
  + \iO (\pi(\phin) + a - b \dt\wn) v
  \non
  \\
  && \quad \hbox{for every $v\in\Vn$ and \aet}\,,
  \label{secondan}
  \\
  \separa
  && \iO \dtt\wn \, v 
  + \iO \nabla( \kuno \dt\wn + \kdue\wn) \cdot \nabla v
  + \lambda \iO \dt\phin \, v
  = \iO g v
  \non
  \\
  && \quad \hbox{for every $v\in\Vn$ and \aet}\,,
  \label{terzan}
  \\
  && \iO \phin(0) \, v 
  = \iO \phiz \, v ,
  \quad 
  \iO \wn(0) \, v
  = \iO \wz \, v ,
  \aand
  \iO \dt\wn(0) \, v
  = \iO \wu \, v,
  \non
  \\
  && \quad \hbox{for every $v\in\Vn$}.
  \label{cauchyn}
\Esist
\Accorpa\Pbln priman cauchyn
\rev{Here, for an arbitrary function $v : Q \to \erre$, we employ an abuse of notation by writing $v(0)$ in place of $v{\vert_{t=0}}(\cdot)$ and we understood it as a function of $x \in \Omega$, that is, $v(0): x \mapsto v(x,0)$. This will be repeatedly used from now  on to simplify the presentation.}
\an{The strategy of the proof can be schematized as follows.}
First, we show that the above problem has a unique solution.
Then, we perform a number of a~priori estimates that allow \an{us} to pass to the limit as $n$ tends to infinity.
In this way, we \an{\juerg{identify} a limit} triple $\soluzeps$\juerg{, which then is shown 
to be a solution to the} problem \Pbleps\ enjoying the desired regularity properties.

\step
Solution to the discrete problem

We represent the unknowns in terms of the basis of the space~$\Vn$.
Namely, we have \aat\ that
\Beq
  \phin(t) = \somma j1n \phinj(t) \ej\,, \quad
  \mun(t) = \somma j1n \munj(t) \ej \,,
  \aand
  \wn(t) = \somma j1n \wnj(t) \ej \,,
  \non
\Eeq
for some functions $\phinj\in H^1(0,T)$, $\munj\in L^2(0,T)$ and $\wnj\in H^2(0,T)$.
Moreover, we introduce the $\erre^n$-valued functions defined \aet~by
\Beq
  \hphin := (\phinj)_{j=1}^n \,, \quad
  \hmun := (\munj)_{j=1}^\an{n} \,,
  \aand
  \hwn := (\wnj)_{j=1}^n \,.
  \non
\Eeq
\juerg{In terms of these true unknowns} the equations \accorpa{priman}{terzan} take the~form
\Bsist
  && \hphin'
  + \an{\AA} \hmun
  + \gamma \hphin
  = \an{\hat {\bf f}}\,,
  \label{primahn}
  \\[2mm]
  && \hmun
  = \an{\AA} \hphin
  + \an{\boldsymbol{\calF}}_\eps(\hphin)
  - b \hwn'\,,
  \label{secondahn}
  \\[2mm]
  && \hwn'' 
  + \an{\AA} (\kuno\hwn'+\kdue\hwn)
  + \lambda \hphin'
  = \an{\hat {\bf g}}\,,
  \label{terzahn}
\Esist
where the matrix $\an{\AA}=(A_{ij})_{i,j=1}^n$
and the vectors $\an{\hat\ff}=(f_i)_{i=1}^n$ and $\an{\hat\gg}=(g_i)_{i=1}^n$ are given~by
\Beq
  A_{ij} := \iO \nabla\ej \cdot \nabla\ei \,, \quad
  f_i := \iO f \ei\,,
  \aand
  g_i := \iO g \ei\,,
  \quad \hbox{for $i,j=1,\dots,n$},
  \non
\Eeq
\gianni{while $\boldsymbol{\calF}_\eps:\erre^n\to\erre^n$ is the function
whose $i$-th component ($i=1,\dots,n$) \juerg{is given by}
\Beq
  \erre^n \ni \rr=(r_1,\dots,r_n) \mapsto
  \iO (\betaeps+\pi) \Bigl( \somma j1n r_j \ej \Bigr) \ei  \,
  \an{+a\iO \ei}
  \,.
  \non
\Eeq
}%
Clearly, $\an{\hat\ff}$~and $\an{\hat\gg}$ are $L^2$ functions
and $\an{\boldsymbol{\calF}}_\eps$~is \Lip\ continuous.
Moreover, the initial conditions \eqref{cauchyn} 
provide initial conditions for \an{the vectors} $\hphin$, $\hwn$ and~$\hwn'$.
We \juerg{first} eliminate $\hphin'$ \juerg{from} \eqref{terzahn} by \an{exploiting} \eqref{primahn}
and then \juerg{eliminate} every occurrence of $\hmun$ by means of~\eqref{secondahn}.
In this way, we obtain a well-posed Cauchy problem for the pair $(\hphin,\hwn)$
coupled with \an{the chemical potential} equation \eqref{secondahn},
and it is clear that the new problem is equivalent to the previous one.
Hence, we find a unique solution with the regularity 
\Beq
  \hphin \in \H1{\erre^n} , \quad
  \hwn \in \H2{\erre^n},
  \aand
  \hmun \in \L2{\erre^n},
  \non
\Eeq
so that the discrete problem has a unique solution, as \an{claimed}.

\medskip

Before \juerg{we start} estimating, we remark a consequence of \an{the compatibility assumptions in}~\eqref{compatibility}.
We choose \juerg{some} $\delta_0>0$ such that both \juerg{the quantities}
$-\rho-(\phizbar)^--\delta_0$ and $\rho+(\phizbar)^++\delta_0$
belong to the interior of~$D(\beta)$.
Then, for some $C_0>0$, we have the inequality
\Bsist
  && \betaeps(r) (r-r_0)
  \geq \delta_0 |\betaeps(r)| - C_0
  \non
  \\
  && \quad \hbox{for every $r\in\erre$, $r_0\in[-\rho-(\phizbar)^-,\rho+(\phizbar)^+]$ and $\eps\in(0,1)$}.
  \qquad
  \label{trickMZ}
\Esist
This is a generalization of \cite[Appendix, Prop.~A.1]{MiZe}.
The detailed proof given in \cite[p.~908]{GiMiSchi} with a fixed $r_0$ 
also works in the present case with \juerg{only} minor changes.

\juerg{Our first estimate} prepares the way to apply the above inequality.

\step
A preliminary estimate

We recall that $\Vn\supset V_1$ and that $V_1$ is the subspace \juerg{of} constant functions.
Hence, we can test \eqref{priman} by $1/|\Omega|$ \an{to} obtain~that
\Beq
  {\phinbar\,}'(t) + \gamma \, \phinbar(t) = \fbar(t)
  \quad \aat,
  \label{primanbar}
\Eeq
whence immediately
\Beq
  \phinbar(t) 
  = \phizbar \, e^{-\gamma t}
  + \iot e^{-\gamma(t-s)} \fbar(s) \, ds
  \quad \hbox{for every $t\in[0,T]$},
  \non
\Eeq
and a simple calculation shows that \pier{(cf.~\eqref{defrho})}
\Beq
  -\rho - (\phizbar)^-
  \leq \phinbar(t)
  \leq \rho + (\phizbar)^+
  \quad \hbox{for every $t\in[0,T]$}.
  \label{stimaphinbar}
\Eeq

Before continuing, it is worth \juerg{making} \an{some} observation\an{s} on projections \juerg{which are} \an{collected in the following remark}.

\Brem
\label{Projection}
Let $\Pn:H\to\Vn$ \an{be} the $H$-\juerg{orthogonal} projection operator.
We list some inequalities that hold true for every $n\in\enne$,
as well as convergence properties as $n$ tends to infinity.
For every $v\in H$, we clearly have that
\Beq
  \norma{\Pn v} \leq \norma v,
  \aand 
  \Pn v \to v 
  \quad \hbox{strongly in $H$}.
  \non
\Eeq
Assume now that $v\in V$.
Then it \an{is} easy to see that \juerg{also}
\Beq
  \Pn v \in V , \quad
  \norma{\nabla\Pn v} \leq \norma{\nabla v},
  \aand
  \normaV{\Pn v} \leq \normaV v \,.
  \non
\Eeq
For a detailed proof, see, e.g., \cite[Rem.~4.2]{CGRS4}.
In particular, we deduce that
\Beq
  \Pn v \to v
  \quad \hbox{strongly in $V$,\quad for every $v\in V$}.
  \non
\Eeq
Next, assume that $v\in W$.
Then, we have~that
\Beq
  v = \somma j1\infty (v,\ej) \ej
  \aand
  - \Delta v = \somma j1\infty (v,\ej) \lambdaj\ej \,.
  \non
\Eeq
We deduce that $\Delta\Pn v=\Pn\Delta v$, and 
we can apply the above inequalities and convergence properties to $\Delta v$ as well
in order to \juerg{recover further information} on~$\Pn v$.
We obtain, with a constant $\CO$ that depends only on~$\Omega$,~that
\begin{alignat*}{2}
  & \norma{\Pn v}_{\Hx2}\leq \CO \norma v_{\Hx2}
  \quad && \hbox{for every $v\in W$}\,,
  \non
  \\[1mm]
  & \norma{\Pn v}_{\Hx3}\leq \CO \norma v_{\Hx3}
  \quad && \hbox{for every $v\in \an{\Hx3\cap W}$}\,,
  \non
  \\[1mm]
  & \Pn v \to v
  \quad && \hbox{strongly in $\Hx2$\quad for every $v\in W$}\,,
  \non
  \\[1mm]
  & \Pn v \to v
  \quad && \hbox{strongly in $\Hx3$\quad for every $v\in \an{\Hx3\cap W}$}.
  \non
\end{alignat*}
Notice that all this can be applied to the initial values
of the discrete solution \an{as they} are projections on~$\Vn$\an{.}
Now, we consider time-dependent functions.
A~simple combination of the above properties with the Lebesgue dominated convergence theorem
shows the following: if we assume that $v\in\L2H$ or $v\in\L2V$ and define $\vn$ by setting
$\vn(t):=\Pn(v(t))$ \aat, then 
\Beq
  \vn \to v
  \quad \hbox{strongly in $\L2H$ or $\L2V$, respectively}.
  \non
\Eeq
\Erem

At this point, we can start estimating,
and we recall that the symbol $c$ stands for possibly different constants independent of $\eps$ and~$n$
\an{a}ccording to our general rule regarding constants stated at the end of Section~\ref{STATEMENT}\an{.}
We repeatedly owe to the properties \PropN\ related to the operator~$\calN$ \an{without further reference}.

\step
First uniform estimate

We \juerg{first} observe that $\calN v\in\Vn$ for every $v\in\Vn$ satisfying $\vbar=0$.
Indeed, both $v$ and $w:=\calN v$ can be expressed in terms of the eigenfunctions~$e_j$, 
and we have that 
\Beq
  \somma {\an{j}}1\infty \lambda_j (w,e_j) \, e_j 
  = -\Delta w
  = v
  = \somma j2n (v,e_j) \, e_j \,.
  \non
\Eeq
Hence, $(w,e_j)=0$ for every $j>n$ (since $\lambda_j>0$ for $j>1$), i.e.,  $w\in\Vn$.
Once this is established, 
we take the difference between \eqref{priman} written for a generic $v\in\Vn$ and \eqref{primanbar} multiplied by~$\iO v$,
write the resulting equality at the time $s\in(0,T)$ and test it by $\calN(\phin-\gianni\phinbar)(s)$.
Then, by integrating \an{over $(0,t)\subset(0,T)$}, we obtain~that
\Bsist
  && \frac 12 \, \normaVp{\phin(t)-\phinbar(t)}^2
  \an{{}+ \intQt \mun (\phin-\phinbar)}
  + \gamma \iot \normaVp{\phin(s)-\phinbar(s)}^2 \, ds
  \non
  \\
  && = \frac 12 \, \normaVp{\phin(0)-\phinbar(0)}^2
  + \intQt (f-\fbar) \calN(\phin-\phinbar).
  \label{testpriman}
\Esist
At the same time, we test \eqref{secondan}, written at the time~$s$, by $-(\phin(s)-\phinbar(s))$
and integrate over~$(0,t)$.
It results~that
\Bsist
  && \intQt |\nabla\phin|^2
  + \intQt \betaeps(\phin) (\phin-\phinbar)
  \an{{}- \intQt \mun (\phin-\phinbar)}
  \non
  \\
  && = - \intQt \bigl( a + \pi(\phin) \bigl) (\phin-\phinbar)
  + b \intQt \dt\wn (\phin-\phinbar) .
  \qquad
  \label{testsecondan}
\Esist
Finally, we take the convolution between \eqref{terzan} and~$1$
(see \eqref{convoluz}) 
and test the resulting equality by~$\dt\wn$.
After time integration, we obtain~that
\Bsist
  && \intQt (\dt\wn-\dt\wn(0)) \dt\wn
  + \frac \kuno 2 \iO |\nabla(\wn(t) - \wn(0))|^2 
  \non
  \\
  && = - \kdue \intQt \nabla(1*\wn) \cdot \nabla\dt\wn
  - \lambda \intQt (\phin-\phin(0)) \dt\wn
  + \intQt (1*g) \dt\wn \,.
  \qquad
  \label{testterzan}
\Esist
At this point, we add \accorpa{testpriman}{testterzan} to each other
and notice that a cancellation occurs.
After rearranging, we deduce~that
\Bsist
  && \frac 12 \, \normaVp{\phin(t)-\phinbar(t)}^2
  + \gamma \iot \normaVp{\phin(s)-\phinbar(s)}^2 \, ds
  \non
  + \intQt |\nabla\phin|^2
  + \intQt \betaeps(\phin) (\phin-\phinbar)
  \non
  \\
  &&  {}
  + \intQt |\dt\wn|^2
  + \frac \kuno 2 \iO |\nabla(\wn(t) - \wn(0))|^2 
  \non
  \\
  \separa
  && = \frac 12 \, \normaVp{\phin(0)-\phinbar(0)}^2
  + \intQt (f-\fbar) \calN(\phin-\phinbar)
  \non
  \\
  &&  {}
  - \intQt \bigl( \pi(\phin) - \pi(\phinbar) \bigl) (\phin-\phinbar)
  - \intQt \bigl( a + \pi(\phinbar) \bigl) (\phin-\phinbar)
  \non 
  \\
  &&  {}
  + b \intQt \dt\wn (\phin-\phinbar) 
  \non 
  + \intQt \dt\wn(0) \dt\wn
  - \kdue \intQt \nabla(1*\wn) \cdot \nabla\dt\wn
  \non 
  \\
  &&  {}
  - \lambda \intQt (\phin-\phinbar) \dt\wn
  - \lambda \intQt (\phinbar-\phin(0)) \dt\wn
  + \intQt (1*g) \dt\wn \,
  \an{=:\sum_{i=1}^{10} I_i}
  \,. 
  \qquad
  \label{sommatestn1}
\Esist
The integral involving $\betaeps$ can be \juerg{estimated} from below 
by combining \an{\eqref{trickMZ} and \eqref{stimaphinbar}} \juerg{as follows:}
\Beq
  \intQt \betaeps(\phin) (\phin-\phinbar)
  \geq \delta_0 \intQt \juerg{|\betaeps(\phin)|}\, - c \,.
  \non
\Eeq
All of the other terms on the \lhs\ are nonnegative.
For those on the \rhs, we \an{perform} separate estimates.

Since the embedding $H\emb\Vp$ is continuous,
the first term \an{$I_1$ is uniformly bounded by the 
assumption~\eqref{hpz} on $\phiz$, Remark~\ref{Projection}, and estimate \eqref{stimaphinbar}.}
Next, we have~that
\begin{align}
  I_2 &= \intQt (f-\fbar) \calN(\phin-\phinbar)
  \leq c \, \norma f_{\L2\Vp}^2 + c \iot \normaV{\calN(\phin-\phinbar)(s)}^2 \, ds
  \non
  \\
  & \leq c \iot \normaVp{(\phin-\phinbar)(s)}^2 \, ds
  + c \,.
  \non
\end{align}
\an{Owing to Young's inequality, the \Lip\ continuity of~$\pi$, and \eqref{stimaphinbar} once more,
we have, for every $\delta>0$,
\begin{align*}
   & I_3 +I_4 +I_5+I_8 
 \leq \frac 14 \intQt |\dt\wn|^2
  + \an{c} \intQt |\phin - \phinbar|^2
  + c
  \\& 
    \leq \frac 14 \intQt |\dt\wn|^2
	+ \delta \intQt |\nabla\phin|^2
  + c_\delta \iot \normaVp{\phin(s)-\phinbar(s)}^2 \, ds \,,
  \non
\end{align*}
where in the second line we also \juerg{used} the compactness inequality~\eqref{compact}.}
\an{Next, \juerg{arguing similarly}, we obtain that
\Bsist
  && I_6 + I_9 + I_{10}
  \leq \frac 14 \intQt |\dt\wn|^2
  + c \intQt \bigl(
    |\dt\wn(0)|^2
    + |\phinbar - \phin(0)|^2
    + |1*g|^2
  \bigr)
  \non
  \\
  && \leq \frac 14 \intQt |\dt\wn|^2
  + c \,,
  \non
\Esist
thanks} to \eqref{stimaphinbar} and to our assumptions on the initial data $\wu$ and~$\phiz$
(by~applying Remark~\ref{Projection}) and on~$g$.
The last term  to be estimated is first treated by an integration by parts.
\gianni{Finally, by} also using Young's inequality and the estimate for $\nabla\wn(0)$ obtained by applying Remark~\ref{Projection}, 
we have, for every $\delta>0$,~that
\Bsist
  && \an{I_7=} - \kdue \intQt \nabla(1*\wn) \cdot \nabla\dt\wn
  \non
  \\
  && = \kdue \intQt |\nabla\wn|^2
  - \kdue \iO \nabla(1*\wn)(t) \cdot \nabla\wn(t)
  \non
  \\
  && \leq c \intQt |\nabla(\wn-\wn(0))|^2
  + \delta \iO |\nabla(\wn(t)-\wn(0))|^2 
  + c_\delta \iO |\nabla(1*\wn)(t)|^2 
  + c \,.
  \non
\Esist
On the other hand, we also have that
\Beq
  \iO |\nabla(1*\wn)(t)|^2 
  = \iO \Bigl| \iot \nabla\wn(s) \, ds \Bigr|^{\pier{2}}
  \leq c \intQt |\nabla\wn|^2
  \leq c \intQt |\nabla(\wn-\wn(0)\an{)}|^2
  + c \,.
  \non
\Eeq
At this point, we recall \eqref{sommatestn1} and all the above estimates,
choose $\delta$ small enough, and apply Gronwall's lemma.
We obtain~that
\Bsist
  && \norma{\phin-\phinbar}_{\L\infty\Vp\cap\L2V}
  + \norma{\betaeps(\phin)}_{\LQ1}
  \non
  \\
  &&  {}
  + \norma{\dt\wn}_{\L2H}
  + \norma{\nabla(\wn-\wn(0))}_{\L\infty H}
  \leq c\,,
  \non
\Esist
whence immediately
\Beq
  \norma\phin_{\L\infty\Vp\cap\L2V}
  + \norma{\betaeps(\phin)}_{\LQ1}
  + \norma\wn_{\H1H\cap\L\infty V}
  \leq c \,.
  \label{primastima}
\Eeq

\step
Consequence

By testing \eqref{secondan} by $1/|\Omega|$, and owing to \eqref{primastima},
we \pier{infer}~that
\Beq
  \norma\munbar_{L^1(0,T)} \leq c \,.
  \label{daprimastima}
\Eeq

\step
Second uniform estimate

We test the equations \eqref{priman}, \eqref{secondan}, and \eqref{terzan}, by
$\mun$, $-\dt\phin$, and $(b/\lambda)\dt\wn$, respectively,
sum up and notice \an{that} the terms involving the products \pier{$\dt\phin\,\mun$} and $\dt\wn\dt\phin$ cancel each other.
Then, we integrate in time \an{and} rearrange \an{to} obtain \an{that}
\Bsist
  && \intQt |\nabla\mun|^2
  + \frac 12 \iO |\nabla\phin(t)|^2
  + \iO \Betaeps(\phin(t))
  \non
  \\
  &&  {}
  + \frac b{2\lambda} \iO |\dt\wn(t)|^2
  + \frac {\kuno b} \lambda \intQt |\nabla\dt\wn|^2
  + \frac {\kuno b} {2\lambda} \iO |\nabla\wn(t)|^2
  \non
  \\
  && = \frac 12 \iO |\nabla\phin(0)|^2
  + \iO \Betaeps(\phin(0))
  + \frac b{2\lambda} \iO |\dt\wn(0)|^2
  + \frac {\kuno b} {2\lambda} \iO |\nabla\wn(0)|^2
  \non
  \\
  && - \gamma \intQt \phin \mun
  + \intQt f \mun
  - \pier{\iO}\Pi(\phin(t)) 
  + \pier{\iO} \Pi(\phin(0)) 
  \non
  \\
  &&  {}
  - a \iO (\phin(t) - \phin(0))
  + \an{\frac b \lambda}\intQt g \dt\wn \,\an{,}
  \label{sommatestn2}
\Esist
\an{where all} \juerg{of the} terms on the \lhs\ are nonnegative.
Moreover, as before, we can recall Remark~\ref{Projection}
in order to estimate the terms involving the initial data,
and just the one containing $\Betaeps$ needs \an{further comments}.
Since $\phiz$ belongs to~$W$ by~\eqref{hpz}, $\phin(0)$~converges to $\phiz$ strongly in~$W$,
\juerg{hence} uniformly.
On the other hand, by the quoted assumption, $\min\phiz$ and $\juerg{{}\max{}}\phiz$ belong to the interior of~$D(\beta)$.
Thus, for some $n_0$ and every $n\geq n_0$, all \juerg{of the} values of $\phin(0)$ 
belong to a compact interval $I$ contained in the interior of~$D(\beta)$.
By also recalling \eqref{disugBetaeps}, we thus may conclude~that
\Beq
  \iO \Betaeps(\phin(0))
  \leq \iO \Beta(\phin(0))
  \leq \max_{r\in I} \Beta(r)
  = c \,.
  \non
\Eeq
It is understood that $n\geq n_0$ from now on, \juerg{which is no restriction since we aim at letting $n$ tend to infinity eventually}.
Let us come to the other terms on the \rhs.
The last one can be dealt with \juerg{employing} Young's inequality
(and then Gronwall's lemma),
and the integral that precedes it has \juerg{already been} estimated,
since it is a multiple of the mean value.
Moreover, since $\Pi$ grows at most quadratically by \an{condition} \eqref{hpPi},
\juerg{we can infer from the compactness inequality \eqref{compact} and \eqref{primastima} that}
\Bsist
  && \iO \Pi(\phin(t)) 
  \leq c \iO |\phin(t)|^2 + c
  \non
  \\
  && \leq \frac 14 \iO |\nabla\phin(t)|^2 
  + c \, \normaVp{\phin(t)}^2
  \an{{}+ c{}}
  \leq \frac 14 \iO |\nabla\phin(t)|^2 
  + c \,.
  \non
\Esist
The other \juerg{integrals} that need some treatment are those containing~$\mun$.
We have~that
\Bsist
  && \intQt (f-\gamma\phin) \mun
  = \intQt (f-\gamma\phin) (\mun-\munbar)
  + \intQt (f-\gamma\phin) \munbar
  \non
  \\
  && \leq \norma{f-\gamma\phin}_{\pier{L^2(0,t;H)}} \, \norma{\mun-\munbar}_{\pier{L^2(0,t;H)}}
  + \norma{f-\gamma\phi}_{\L\infty\Vp} \, \norma\munbar_{\L1V}
  \non
  \\[1mm]
  && \leq c \, \norma{\nabla\mun}_{\pier{L^2(0,t;H)}}
  + c \, \norma\munbar_{L^1(0,T)}
  \leq \pier{{}\frac12 \intQt |\nabla\mun|^2 + c{}}\,,
  \non
\Esist
where we have used the Poincar\'e inequality \eqref{poincare},
our assumptions on~$f$ (see~\eqref{hpfg}), \eqref{primastima}, and~\eqref{daprimastima}.
By coming back to \eqref{sommatestn2}, collecting the above estimates and observations, 
and applying the Gronwall lemma, 
we conclude~that
\Bsist
  && \norma{\nabla\mun}_{\L2H}
  + \norma\phin_{\L\infty V}
  + \norma{\Betaeps(\phin)}_{\L\infty\Luno}
  \non
  \\
  &&  {}
  + \norma\wn_{\W{1,\infty}H\cap\H1V}
  \leq c \,.
  \label{secondastima}
\Esist

\step
Third uniform estimate

\juerg{Next,} recalling that every constant is allowed as a test function, 
we test \eqref{secondan} by $\phi\an{_n}(t)-\phibar\an{_n}(t)$ and rearrange.
Omitting the time variable for brevity, we have \aet\ that
\Bsist
  && \iO |\nabla\phin|^2
  + \iO \betaeps(\phin)(\phin-\phinbar)
  \non
  \\
  && = \iO \mun (\phin-\phinbar)
  - \iO \pi(\phin) (\phin-\phinbar)
  - a \iO (\phin-\phinbar)
  \non
  \\
  &&  {}
  + b \iO \dt\wn (\phin-\phinbar) .
  \label{testsecondan3}
\Esist
\juerg{In view of} \eqref{stimaphinbar} and of our assumption \eqref{compatibility},
we can bound the integral involving $\betaeps$ from below using \eqref{trickMZ}: 
\Beq
  \iO \betaeps(\phin)(\phin-\phinbar)
  \geq \delta_0 \iO |\betaeps(\phin)| - c \,.
  \non
\Eeq
As for the \rhs, the first term \juerg{needs} some treatment.
Thanks to Poincar\'e's inequality \eqref{poincare} and to \eqref{secondastima}, we have~that
\Beq
  \iO \mun (\phin-\phinbar)
  = \iO (\mun-\munbar) (\phin-\phinbar)
  \leq c \, \norma{\nabla\mun} \, \norma{\phin-\phinbar}
  \leq c \, \norma{\nabla\mun} \,.
  \non
\Eeq
The sum of the other  \juerg{terms} is bounded from above~by
\Beq
  c \bigl(
    \norma\phin^2
    + |\phinbar|^2  
    + \norma{\dt\wn}^2
  \bigr)
  + c \,.
  \non
\Eeq
\juerg{Combining} \eqref{testsecondan3}, the inequalities just derived, and the previous estimates,
we see that \juerg{the function} $\juerg{t\mapsto\iO|\betaeps(\varphi_n(t))}|$ is bounded from above by an $\juerg{L^2(0,T)}$ function independently of both $n$ and~$\eps$,
that is, it holds 
\Beq
  \norma{\betaeps(\phin)}_{\L2\Luno} \leq c \,,
  \label{perterzastima}
\Eeq
whence we trivially derive an estimate in $L^2(0,T)$ for the mean value of~$\betaeps(\phin)$.
Then, from \eqref{secondan}, we can estimate the \juerg{$L^2(0,T)$ norm of~$\munbar$}.
This, \eqref{secondastima}, and the use of the Poincar\'e inequality once more, imply~that
\Beq
  \norma\mun_{\L2V} \leq c \,.
  \label{terzastima}
\Eeq

\step
\juerg{Fourth} uniform estimate

We recall Remark~\ref{Projection} and use \an{the} notations \an{introduced there}.
We \an{fix} some $v\in\L2V$ and define $\vn\in\L2\Vn$ by setting 
$\vn(t)=\Pn(v(t))$ \aat.
Then, we test \eqref{priman} by $\vn$, and integrate \an{over} time \an{to} obtain~that
\Beq
  \intQ \dt\phin \, \vn
  = - \intQ \nabla\mun \cdot \nabla\vn
  + \intQ (f - \gamma\phin) \vn
  \leq c \norma\vn_{\L2V}
  \leq c \norma v_{\L2V}\an{.}
  \non
\Eeq
On the other hand, we have that
\Beq
  \intQ \dt\phin \, \vn
  = \intQ \dt\phin \, v\,,
  \non
\Eeq
since $\dt\phin$ is $\Vn$-valued.
\an{Thus,} we \an{readily} conclude that
\Beq
  \norma{\dt\phin}_{\L2\Vp} \leq c \,.
  \label{quartastima}
\Eeq

\step
Fifth uniform estimate
 
\pier{Thanks to \eqref{quartastima}, the} same argument, applied to 
\an{equation}~\eqref{terzan} for $\wn$, yields that
\Beq
  \norma{\dtt\wn}_{\L2\Vp} \leq c \,.
  \label{quintastima}
\Eeq

\step
\an{\juerg{Passage} to the l}imit 

At this point, we can pass to the limit as $n\to\infty$.
Indeed, by recalling \eqref{primastima}, \eqref{secondastima} and \accorpa{terzastima}{quintastima},
and applying well-known weak, weak star, and strong compactness results
(for the latter see, e.g., \cite[Sect.~8, Cor.~4]{Simon}),
we \pier{deduce} that there exists a triple $\soluzeps$ such~that
\Bsist
  & \phin \to \phieps
  & \quad \hbox{weakly star in $\H1\Vp\cap\L\infty V$}
  \non
  \\
  && \qquad \hbox{and strongly in $\C0H$}\,,
  \label{convphin}
  \\
  & \mun \to \mueps
  & \quad \hbox{weakly in $\L2V$}\,,
  \label{convmun}
  \\
  & \wn \to \weps 
  & \quad \hbox{weakly star in $\H2\Vp\cap\W{1,\infty}H\cap\H1V$}
  \non
  \\
  && \qquad \hbox{and strongly in $\pier{{}\H1H{}}\cap\C1\Vp$}\,, 
  \label{convwn}
\Esist
as $n$ tends to infinity (at least for a subsequence which is not relabeled).
Moreover, since $\betaeps$ and $\pi$ are \Lip\ continuous,
we also have~that
\Beq
  \betaeps(\phin) \to \betaeps(\phieps)
  \aand 
  \pi(\phin) \to \pi(\phieps)
  \quad \hbox{strongly in $\C0H$}.
  \label{convbetaepsphin}
\Eeq
We claim that \juerg{this}  triple is a (weak) solution to problem \Pbleps.
Since $\phin(0)$, $\wn(0)$ and $\dt\wn(0)$ are the $H$ projections of $\phiz$, $\wz$ and~$\wu$, 
they strongly converge in $H$ to $\phiz$, $\wz$ and~$\wu$, respectively.
On the other hand, they converge to $\phieps(0)$, $\weps(0)$ and~$\dt\weps(0)$, respectively, 
strongly (at~least) in~$\Vp$, thanks to \eqref{convphin} and~\eqref{convwn}.
Hence, the initial conditions \eqref{cauchyeps} are satisfied.
Now, we show that the variational equations \accorpa{primaeps}{terzaeps} are satisfied as well.
We recall Remark~\ref{Projection}\an{, fix} any $v\in\L2V$, define $\vn\in\L2\Vn$ by setting $\vn(t):=\Pn(v(t))$ \aat,
and observe that $\vn$ converges to $v$ strongly in $\L2V$.
Next, we test each of the equations \accorpa{priman}{terzan} by~$\vn$ and integrate in time over $(0,T)$.
At this point, on account of the convergence properties proved or mentioned,
it is \sfw\ to pass to the limit as $n\to\infty$ in the equalities we obtain.
The resulting equalities are the same equations with $\soluzeps$ in place of $\soluzn$, i.e.,
the time-integrated versions of \accorpa{priman}{terzan} with arbitrary time-dependent test functions $v\in\L2V$,
which are equivalent to \accorpa{priman}{terzan} themselves.

\step
Conclusion of the proof

It remains to establish the stronger regularity requirements \an{stated in \eqref{piureg}}. 
To this end, we see that, \aet, $\phieps(t)$ is a solution $u\in V$ to the nonlinear elliptic problem
\Beq
  \iO \nabla u \cdot \nabla v
  + \iO \betaeps(u) v
  = \iO h v
  \quad \quad \hbox{for every $v\in V$}\,,
  \label{elliptic}
\Eeq
where $h$ \an{is} the value of $\mueps-\pi(\phieps)-a+b\dt\weps$ \an{evaluated} at~$t$ in our case.
On the other hand, every solution $u$ to problem \eqref{elliptic} 
satisfies the estimate
\Beq
  \norma{\betaeps(u)}_6 \leq \norma h_6 \,,
  \label{stimasei}
\Eeq
whenever $h\in\Lx6$.
To show that \eqref{stimasei} actually holds true,
we can formally choose $v=(\betaeps(u))^5$ in \eqref{elliptic}
(to~be more rigorous, we should use \an{a suitable} truncation).
Next, we apply the generalized Young inequality 
with conjugate exponents $6$ and $6/5$ to the resulting \rhs\ and rearrange.
Then, \eqref{stimasei} plainly follows.
Moreover, by elliptic regularity, we infer~that
\Beq
  u \in \an{\Wx{2,6}\cap W} 
  \aand
  \norma u_{\Wx{2,6}}
  \leq \CO \bigl( \norma u + \norma h_6)\,,
  \label{regelliptic}
\Eeq
with a constant $\CO$ that depends only on~$\Omega$.
\an{Then we apply \eqref{stimasei} with $u=\phieps(t)$, square} and integrat\an{e} in time \an{to} deduce that
\Beq
  \norma{\betaeps(\phieps)}_{\L2{\Lx6}}^2
  \leq \norma{\mueps-\pi(\phieps)-a+b\dt\weps}_{\L2{\Lx6}}^2 \,.
  \non
\Eeq
Since the \rhs\ of this inequality is uniformly bounded \an{owing} to our previous estimates
and the continuous embedding $V\emb\Lx6$,
we conclude~that
\Beq
  \betaeps(\phieps) \in \L2{\Lx6}
  \aand
  \norma{\betaeps(\phieps)}_{\L2{\Lx6}} \leq c \,.
  \label{stimabetaeps}
\Eeq
Similarly, by applying \eqref{regelliptic}, we also have~that
\Beq
  \phieps \in \L2{\Wx{2,6}}
  \aand
  \norma\phieps_{\L2{\Wx{2,6}}} \leq c \,.
  \label{stimaW26}
\Eeq


\section{Existence and regularity}
\label{EXISTENCE}
\setcounter{equation}{0}
\an{This final part of the paper is \gianni{devoted} 
to prove the \pier{existence and} regularity results stated in} Theorems~\ref{Existence} and~\ref{Regularity}.

\subsection{Proof of Theorem~\ref{Existence}}

\an{To proceed rigorously, let us consider the discrete problem \Pbln\ analyzed in the previous section.}
By the lower semicontinuity of \an{norms},
it is clear that the bounds 
\eqref{primastima}, \eqref{secondastima}, and \accorpa{terzastima}{quintastima},
proved for the discrete solution $\soluzn$
are conserved with the same constants in the limit as $n\to\infty$.
By also accounting for \accorpa{stimabetaeps}{stimaW26},
we thus have~that
\Bsist
  && \norma\phieps_{\H1\Vp\cap\L\infty V\cap\L2{\Wx{2,6}}}
  + \norma\mueps_{\L2V}
  \non
  \\
  &&
  + \norma{\betaeps(\phieps)}_{\L2{\Lx6}}
  + \norma\weps_{\H2\Vp\cap\W{1,\infty}H\cap\H1V}
  \leq c\,,
  \label{stimaeps}
\Esist
and we recall that, according to our general rule, 
the constant $c$ in the \an{above line} has the same dependence 
as the constant $K_1$ of the statement\juerg{. In particular, it is independent of $\eps$}.
From \eqref{stimaeps} and the compactness results already mentioned, we have~that
\Bsist
  & \phieps \to \phi
  & \quad \hbox{weakly star in $\H1\Vp\cap\L\infty V$}
  \non
  \\
  && \qquad \hbox{and strongly in $\C0H$}\,,
  \label{convphieps}
  \\
  & \mueps \to \mu
  & \quad \hbox{weakly in $\L2V$}\,,
  \label{convmueps}
  \\
  & \betaeps(\phieps) \to \xi
  & \quad \hbox{weakly in $\L2{\Lx6}$}\,,
  \label{convxieps}
  \\
  & \weps \to w
  & \quad \hbox{weakly star in $\H2\Vp\cap\W{1,\infty}H\cap\H1V$}
  \non
  \\
  && \qquad \hbox{and strongly in $\pier{{}\H1H{}}\cap\C1\Vp$}\,, 
  \label{convweps}
\Esist
for some quadruplet $\soluz$ as $\eps$ tends to zero 
(at least for a \an{not relabeled subsequence}).
Notice that this quadruplet satisfies the estimate \eqref{stima}
by the lower semicontinuity of \an{norms}.
We now prove that it is a solution to problem \Pbl.
Clearly, the initial conditions \eqref{cauchy} are \an{fulfilled}.
Moreover, \pier{by the maximal monotonicity of $\beta$} \an{it is a standard matter to realize that} the condition $\xi\in\beta(\phi)$ that appears in \eqref{seconda} is satisfied as well.
Indeed, it suffices to combine the strong convergence of $\phieps$,
the weak convergence of $\betaeps(\phieps)$, and a well-known property of the Yosida approximation
(see, e.g.,  \cite[Prop.~2.2, p.~38]{Barbu}).
Finally, as in the previous proof, it is \sfw\ to pass to the limit
in the time-integrated versions of the equations \accorpa{primaeps}{terzaeps}
\juerg{in order to} obtain the time-integrated versions of the equations \accorpa{prima}{terza}
with arbitrary time-dependent test functions $v\in\L2V$,
which are equivalent to \accorpa{prima}{terza} themselves.
This completes the proof.

\subsection{Proof of Theorem~\ref{Regularity}}

Following the line \an{of arguments} of the proof of Theorem~\ref{Existence},
we use the estimates already established 
for the discrete solution $\soluzn$ and the approximating solution $\soluzeps$
and perform further estimates. 
So, we want to show~that
\Beq
  \norma\phin_{\H1V}
  + \norma\mun_{\L\infty V}
  + \norma\wn_{\pier{\H2H\cap\W{1,\infty}V}}
  \leq c\,, 
\label{stimaregn}
\Eeq
at least for sufficiently large $n\in\enne$, as well~as
\Beq
  \norma{\betaeps(\phieps)}_{\L\infty{\Lx6}} \leq c
  \aand
  \norma\phieps_{\L\infty{\Wx{2,6}}} \leq c \,,
  \label{stimaregeps}
\Eeq
with a constant $c$ that has the same dependence as the constant $K_3$ in the statement.
To prove \eqref{stimaregn}, we first observe that the component $\mun$ of the discrete solution $\soluzn$
is more regular than required:
\juerg{indeed, it is \Lip\ continuous, 
as follows from looking at $\hmun$} in equation \eqref{secondahn}.
Hence, we are allowed to take $t=0$ in \eqref{secondan}.
It results~that
\Beq
  \mun(0) = \Pn(-\Delta\phiz+\betaeps(\phin(0))+\pi(\phin(0))+a-b\wu)\,,
  \label{defmuz}
\Eeq
where $\Pn:H\to\Vn$ is the orthogonal projection operator.
Recall that $\mun(0)$ also depends on~$\eps$, of course, despite of the \an{used notation}.
It is convenient to first establish an estimate for~$\mun(0)$. 

\Blem
\label{Stimamuz}
There exist a positive constant $c$
and a positive integer $n_0$ such that the inequality
\Beq
  \normaV{\mun(0)} \leq c
  \label{stimamuz}
\Eeq
holds true for every $\eps\in(0,1)$ and every $n\geq n_0$.
\Elem

\Bdim
First, we prove that
\Beq
  \normaV{\betaeps(\phin(0))} \leq c 
  \label{technical}
\Eeq
for every $\eps\in(0,1)$, some $n_0$ and and every $n\geq n_0$.
By recalling \eqref{hpz} \pier{and \eqref{compatibility}}, we can find elements $\rmin$ and $\rmax$ \juerg{in} the interior of $D(\beta)$
satisfying $\rmin<0<\rmax$, $\rmin<\juerg{\min}\,\phiz$ and $\rmax>\juerg{\max}\,\phiz$.
Next, since $\phiz\in W$, Remark~\ref{Projection} ensures that
$\phin(0)$ converges to $\phiz$ in $\Hx2$ as $n\to\infty$, thus uniformly.
Therefore, there exists \juerg{some $n_0\in\enne$} such that $\rmin\leq\phin(0)\leq \rmax$ for every $n\geq n_0$,
so~that (by~\eqref{disugbetaeps})
\Beq
  |\betaeps(\phin(0))|
  \leq \sup_{s\in[\rmin,\rmax]} |\beta(s)|
  \quad \hbox{in $\Omega$, \ for every $n\geq n_0$}.
  \non
\Eeq
In particular, \juerg{the sequence $\{\betaeps(\phin(0))\}$} is uniformly bounded in~$H$.
On the other hand, since the restriction of $\beta$ to the interior of $D(\beta)$ is a $C^1$ function by \eqref{hpbetareg},
the following inequality holds:
\Beq
  |\betaeps'(r)| \leq \sup_{s\in[\rmin,\rmax]} |\beta'(s)| =:C
  \quad \hbox{for every $r\in[\rmin,\rmax]$} .
  \non  
\Eeq
For a detailed proof (with a different notation) see, e.g., \cite[formula (5.2)]{CGRS4}.
Then, we have that
\Beq
  \norma{\nabla\betaeps(\phin(0))}
  = \norma{\betaeps'(\phin(0)) \nabla\phin(0)}
  \leq C \, \norma{\nabla\phiz} \,,
  \non
\Eeq
\an{so that \eqref{technical} follows.}
At this point, we easily derive~\eqref{stimamuz}.
By also accounting for assumption~\eqref{hpdatireg} and Remark~\ref{Projection} once more, 
we have indeed
\Bsist
  && \normaV{\mun(0)}
  \leq \normaV{-\Delta\phiz+\betaeps(\phin(0))+\pi(\phin(0))+a-b\wu}
  \non
  \\[1mm]
  && \leq c \bigl( \norma\phiz_{\Hx3} + \normaV{\betaeps(\phin(0))} + \normaV{\phin(0)} + \normaV\wu + 1 \bigr)
  \leq c \,.
  \non
\Esist
\Edim

\an{Let us now} continue \an{with} the proof.
It is understood that $n\geq n_0$ (given by the lemma) from now~on.
In order to make the argument more transparent,
it is convenient to prepare an auxiliary estimate
depending on a positive parameter $M$ whose value will be chosen later~on.

\step
Auxiliary estimate

We repeat part of the argument used to arrive at \eqref{perterzastima},
but \juerg{this time we} avoid time integration.
We account for \eqref{stimaphinbar} in order to apply \eqref{trickMZ} once more.
We test \eqref{secondan} \aet\ by $M(\phin-\phinbar)$.
Then, we \juerg{invoke the Poincar\'e inequality \an{\eqref{poincare}}
and} the Young inequality \eqref{young} with  $\delta=(8M\CO)^{-1}$.
By also taking advantage of~\eqref{secondastima},
we find (\aet)~that
\Bsist
  && \delta_0 M |\Omega| |\overline{\betaeps(\phin)}|
  \leq \delta_0 M \iO |\betaeps(\phin)|
  \leq M \Bigl(\iO \betaeps(\phin) (\phin-\phinbar) + C_0|\Omega| \Bigr)
  \non
  \\
  && \leq M \iO (\mun-\munbar) (\phin-\phinbar) 
  - M \iO \pi(\phin) (\phin-\phinbar)
  \non
  \\
  &&  {}
  - M \iO (a-b\dt\wn) (\phin-\phinbar)
  + c_M
  \non
  \\
  && \leq \frac 18 \iO |\nabla\mun|^2
  + c_M \bigl( \norma\phin^2 + \norma{\dt\wn}^2 + |\phinbar|^2 + 1 \bigr)
  \non
  \\
  && \leq \frac 18 \iO |\nabla\mun|^2 + c_M \,\an{,}
  \non
\Esist
\an{with \pier{the} positive constant $C_0$ arising from \eqref{trickMZ}.}
Since \eqref{secondan} and the \juerg{already known} estimates for $\phin$ and $\dt\wn$ imply that
\begin{align*}
	|\munbar|\leq|\overline{\betaeps(\phin)}|+c 
	\quad  \aet ,
\end{align*}
we deduce~that
\Beq
  \delta_0 M |\Omega| |\munbar(t)|
  \leq \frac 18 \iO |\nabla\mun(t)|^2 + c_M 
  \quad \aat .
  \label{compensation}
\Eeq

\step
Sixth uniform estimate

\an{By virtue of the already proved regularity of~$\mun$,}
we can \an{now} take $\dt\mun$ as \juerg{a} test function in~\eqref{priman} and,
at the same time, we can differentiate \eqref{secondan} with respect to time
and then test the \juerg{resulting equality} by \pier{$- \dt\phin$}.
We do this and also test \eqref{terzan} by $(b/\lambda)\dtt\wn$.
Then, we sum up and notice that four terms cancel each other.
Finally, we integrate with respect to time and add \eqref{compensation} to the \an{resulting} equality.
\an{Collecting the terms\juerg{, we obtain that}}
\Bsist
  && \frac 12 \iO |\nabla\mun(t)|^2
  + \intQt |\nabla\dt\phin|^2
  + \intQt \betaeps'(\phin)|\dt\phin|^2
  \non
  \\
  &&  {}
  + \frac b\lambda \intQt |\dtt\wn|^2
  + \frac {b\kuno}{2\lambda} \iO |\nabla\dt\wn(t)|^2
  + \delta_0 M |\Omega| |\munbar(t)|
  \non
  \\
  && \leq \frac 12 \iO |\nabla\mun(0)|^2
  + \intQt (f-\gamma\phin) \dt\mun
  - \intQt \pi'(\phin)|\dt\phin|^2
  \non
  \\
  &&  {} 
  + \frac {b\kuno}{2\lambda} \iO |\nabla\dt\wn(0)|^2
  - \frac{b\kdue}\lambda \intQt \nabla\wn \cdot \nabla\dtt\wn
  + \frac b\lambda \intQt g \, \dtt\wn 
  \non
  \\
  &&  {}
  + \frac 18 \iO |\nabla\mun(t)|^2
  + c_M \,\an{,}
  \label{regtest}
\Esist
\an{where all} \juerg{of} the terms on the \lhs\ are nonnegative.
The first term on the \rhs\ is estimated by the \an{above} lemma.
The other term involving an initial value is bounded by the $V$-norm of~$\wu$.
As for the first volume term on the \rhs, we integrate by parts in time and have~that
\Bsist
  && \intQt (f-\gamma\phin) \dt\mun
  \non
  \\
  && = - \intQt (\dt f - \gamma\dt\phin) \mun
  + \iO (f-\gamma\phin)(t) \mun(t)
  - \iO (f-\gamma\phin)(0) \mun(0).
  \non 
\Esist
The volume integral \an{on the right} is estimate\an{d} by $\norma{\dt f-\gamma\dt\phin}_{\L2\Vp} \, \norma\mun_{\L2V}$\juerg{,
which is bounded} on account of \eqref{hpdatireg}, \an{\eqref{terzastima} and \eqref{quartastima}}.
The last term is easily treated \an{once again with the help of} the lemma.
The remaining term is dealt with by using the Young and Poincar\'e inequalities:
\Bsist
  && \iO (f-\gamma\phin)(t) \mun(t)
  = \iO (f-\gamma\phin)(t) (\mun-\munbar)(t)
  + \iO (f-\gamma\phin)(t) \munbar(t)
  \non
  \\
  && \leq \frac 18 \iO |\nabla\mun(t)|^2
  + c \, \norma{f-\gamma\phin}_{\L\infty H}^2
  + C^* \, |\munbar(t)|
  \non
  \\
  && \leq \frac 18 \iO |\nabla\mun(t)|^2
  + C^* \, |\munbar(t)|
  + c\,,
  \label{marked}
\Esist
where we have used the special symbol $C^*$ to mark the constant in front of $|\munbar(t)|$ for \an{future} reference.
Notice that $C^*$ is a multiple of an upper bound for the norm of $\norma{f-\gamma\phin}$ in $\L\infty\Vp$,
which is known by \eqref{hpdatireg} and \eqref{primastima}.
Next, \pier{it turns out}~that
\Beq
  - \intQt \pi'(\phin)|\dt\phin|^2
  \leq \frac 12 \intQt |\nabla\dt\phin|^2 
  + c \norma{\dt\phin}_{\L2\Vp}^2
  \leq \frac 12 \intQt |\nabla\dt\phin|^2 
  + c\,,
  \non
\Eeq
thanks to the \Lip\ continuity of~$\pi$, the compactness inequality \eqref{compact}, and \eqref{quartastima}.

It remains to estimate the volume integrals involving $\wn$ that appear \an{on the \rhs\ of}~\eqref{regtest}.
The last one is trivially treated via Young's inequality.
The other can be dealt with \juerg{as follows, using} integration by parts, \eqref{secondastima}, 
the Young inequality, and Remark~\ref{Projection}\pier{. Indeed, we have that}
\Bsist
  && - \frac{b\kdue}\lambda \intQt \nabla\wn \cdot \nabla\dtt\wn 
  \non
  \\
  && = \frac{b\kdue}\lambda \intQt |\nabla\dt\wn|^2
  - \frac{b\kdue}\lambda \iO \nabla\wn(t) \cdot \nabla\dt\wn(t)
  + \frac{b\kdue}\lambda \iO \nabla\wn(0) \cdot \nabla\dt\wn(0) 
  \non
  \\
  && \leq c \, \norma\wn_{\H1V}^2
  + \frac {b\kuno}{4\lambda} \iO |\nabla\dt\wn(t)|^2
  + c \, \norma\wn_{\L\infty V}^2
  + c \, \normaV\wz \, \normaV\wu 
  \non
  \\
  && \leq \frac {b\kuno}{4\lambda} \iO |\nabla\dt\wn(t)|^2
  + c \,.
  \non
\Esist
At this point, we can easily conclude.
Indeed, if we choose $M$ in order that $\delta_0 M|\Omega|=C^*+1$,
and collect \eqref{regtest} and all the above estimates, \juerg{then}
we obtain
\gianni{%
\Beq
  \norma\phin_{\H1V}
  + \norma\mun_{\L\infty V}
  + \norma\wn_{\H2H\cap\W{1,\infty}V}
  \leq c
  \label{sestastima}
\Eeq
}%
with a constant $c$ that has the same dependence 
on the structure \juerg{and} the data as required, since even $M$ has this property.

\step
Conclusion of the proof

\an{We are then left with checking \eqref{stimaregeps}.}
To this end, it suffices to come back to the nonlinear elliptic problem \eqref{elliptic}
and the corresponding estimates \accorpa{stimasei}{regelliptic}
and to argue as we did to prove \accorpa{stimabetaeps}{stimaW26},
in this case avoiding time integration.

At this point, we can easily conclude.
As the discrete solution $\soluzn$ converges as $n\to\infty$ to the solution $\soluzeps$ to the approximating problem,
it is clear that the analogue of \eqref{stimaregn} for $\soluzeps$ holds true with the same constant,
by the semicontinuity of \an{norms}.
We conclude that the convergence properties \accorpa{convphieps}{convweps}
can be improved on account of the estimate just mentioned and \eqref{stimaregeps}.
On the other hand, the previous proof ensures that the limiting quadruplet $\soluz$
is a solution to problem \Pbl, 
and the estimates proved for the approximating solution are conserved in the limit.
Therefore, the proof of Theorem~\ref{Regularity} is complete.


\vskip 6mm
\noindent{\bf Acknowledgments}

\noindent
This research was supported by the Italian Ministry of Education, 
University and Research (MIUR): Dipartimenti di Eccellenza Program (2018--2022) 
-- Dept.~of Mathematics ``F.~Casorati'', University of Pavia. 
In addition, {PC and AS gratefully acknowledge some other support 
from the MIUR-PRIN Grant 2020F3NCPX ``Mathematics for industry 4.0 (Math4I4)'' and}
their affiliation to the GNAMPA (Gruppo Nazionale per l'Analisi Matematica, 
la Probabilit\`a e le loro Applicazioni) of INdAM (Isti\-tuto 
Nazionale di Alta Matematica). 


\End{document}
